\documentclass[12pt]{amsart}
\usepackage{amsfonts}
\usepackage{epsfig}
\usepackage{graphicx}
\usepackage{amsmath}
\usepackage{amssymb}
\usepackage{latexsym}
\usepackage{rotating}
\usepackage{pstricks, pst-node, pst-text, pst-3d}
\usepackage{amsbsy}

\usepackage{amsthm}
\usepackage{mathrsfs}
\usepackage[all]{xy}
\usepackage{a4wide}

\input{amssym.def}

\DeclareMathOperator{\real}{\mathbb{R}}
\DeclareMathOperator{\comp}{\mathbb{C}}

\newtheorem{theorem}{Theorem}
\newtheorem{proposition}{Proposition}
\newtheorem{corollary}{Corollary}
\newtheorem{lemma}{Lemma}
\newtheorem{definition}{Definition}
\theoremstyle{remark}
\newtheorem{remark}{Remark}

\DeclareMathOperator{\re}{Re }
\DeclareMathOperator{\im}{Im }
\DeclareMathOperator{\sgn}{sgn }
\DeclareMathOperator{\Arg}{Arg}
\DeclareMathOperator{\AdS}{AdS}
\DeclareMathOperator{\spn}{span} 
\DeclareMathOperator{\pr}{pr} 

\DeclareMathOperator{\SU}{SU}
\DeclareMathOperator{\U}{U}

\DeclareMathOperator{\wSU}{\widetilde{SU}}
 
\DeclareMathOperator{\SL}{SL}

\DeclareMathOperator{\scrS}{\mathscr S}
\DeclareMathOperator{\scrC}{\mathscr C}
\DeclareMathOperator{\scrT}{\mathscr T}

\input epsf
\begin{document}
\title[Sub-Riemannian and sub-Lorentzian geometry on $\SU(1,1)$...]{Sub-Riemannian and sub-Lorentzian geometry~on~$\SU(1,1)$
and on its universal cover}
\author[E.~Grong and A.~Vasil'ev]{Erlend Grong and Alexander Vasil'ev}

\thanks{The authors have been  supported by the grant of the Norwegian Research Council \#204726/V30, by the NordForsk network `Analysis and Applications', grant \#080151, and by the European Science Foundation Research Networking Programme HCAA}


\address{\newline E.~Grong and A.~Vasil'ev: \smallskip \newline Department of Mathematics \newline University of Bergen \newline  P.O.~Box~7803
\newline Bergen N-5020 \newline Norway \smallskip}
\email{erlend.grong@math.uib.no}
\email{alexander.vasiliev@math.uib.no}

\begin{abstract}
We study sub-Riemannian and sub-Lorentzian geometry on the Lie group $\SU(1,1)$ and on its universal cover $\wSU(1,1)$. In the sub-Riemannian case we find the distance function and completely describe sub-Riemannian geodesics on both $\SU(1,1)$ and  $\wSU(1,1)$,
connecting two fixed points. In particular, we prove that there is a strong connection between the conjugate loci and the number of geodesics. In the sub-Lorentzian case, we describe the geodesics connecting two points on  $\wSU(1,1)$, and compare them with Lorentzian ones. It turns out that the reachable sets for Lorentzian and sub-Lorentzian normal geodesics intersect but are not included one to the other.  A description of
the timelike future is obtained and compared in the Lorentzian and sub-Lorentzain cases. 
\end{abstract}
\maketitle

\vspace{-10pt}


\centerline{\sc Contents}

\qquad 1. Introduction.\smallskip

\qquad 2.  Sub-Riemannian and sub-Lorentzian geometry. 

\qquad \hspace{10pt} {\small 2.1.  Sub-Riemannian manifolds.}

\qquad \hspace{10pt} {\small 2.2.  Sub-Lorentzian manifolds.}

\qquad \hspace{10pt} {\small 2.3.  Minimizing and maximizing curves seen from the viewpoint of optimal control.}

\smallskip

\qquad 3.    Structure of the Lie groups $\SU(1,1)$ and $\wSU(1,1)$.

\qquad \hspace{10pt} {\small 3.1.  Some notations.}

\qquad \hspace{10pt} {\small 3.2.  Lie group $\SU(1,1)$.}

\qquad \hspace{10pt} {\small 3.3.  Lie group structure of the universal cover $\wSU(1,1)$ of $\SU(1,1)$.}

\smallskip

\qquad 4. Sub-Riemannian geometry on $\SU(1,1)$ and $\wSU(1,1)$.

\qquad \hspace{10pt} {\small 4.1.  Geodesics, horizontal space, and vertical space.}

\qquad \hspace{10pt} {\small 4.2.  Length and number of geodesics.}

\qquad \hspace{10pt} {\small 4.3.  The cut and conjugate loci.}

\smallskip

\qquad 5. Sub-Lorentzian geometry on $\wSU(1,1)$.

\qquad \hspace{10pt} {\small 5.1.  Sub-Lorentzian maximizers and geodesics on $\wSU(1,1)$.}

\qquad \hspace{10pt} {\small 5.2.  Number of geodesics.}

\qquad \hspace{10pt} {\small 5.3.  Lorentzian and sub-Lorentzian timelike future.}

\smallskip

\qquad 6.  Proofs of main results.\smallskip

\qquad \hspace{10pt} {\small 6.1.  Proof of Proposition 2.}

\qquad \hspace{10pt} {\small 6.2.  Proof of Corollary 1.}

\qquad \hspace{10pt} {\small 6.3.  Proof of Proposition 5.}

\section{Introduction}

Sub-Riemannian geometry is proved to play an important role in many applications, e.g., in
mathematical physics and control theory.
Sub-Riemannian geometry enjoys major differences from the Riemannian being a 
generalization of the latter at the same time, e.g., geodesics may be singular, the 
Hausdorff dimension
is larger than the manifold topological dimension, the exponential map is never a local 
diffeomorphism. There exists a large amount of literature developing sub-Riemannian 
geometry. Typical general references are \cite{Mon, Str1, Str2}.  The sub-Lorentzian case 
is less studied and the first works in this directions appeared rather recenlty, see
\cite{ChangMarkVas, Groch1, Groch2, Korolko}.

In the development of sub-Riemannian geometry, one observes several  examples, which are mainly nilpotent Lie groups,
with either a left or right invariant distribution and metric. A sample representative is the Heisenberg group (see, e.g., \cite{Mon}). Analysis of these groups in the sub-Riemannian setting,
is already well studied.
While these groups enjoy the advantage that explicit results are easier to obtain,
their properties are sometimes too nice to be good examples to reveal all specific features of
sub-Riemannian geometry in its generality. For instance, the cut locus and the conjugate loci for the Heisenberg group globally
coincide.

A natural next step after considering nilpotent groups is to consider semisimple Lie groups. Let $G$ be such a group with the Lie algebra
$\mathfrak{g}$. Let $\iota: \mathfrak{g} \rightarrow \mathfrak{g}$ be a Cartan involution (ex. for matrix algebras
we can define $\iota$ as the map that sends an element $X$ to minus its conjugate transpose). Then there is a splitting
$\mathfrak{g} = \mathfrak{k} \oplus \mathfrak{p}$, where $\mathfrak{k}$ and $\mathfrak{p}$
are the $+1$ and $-1$ eigenspaces of $\iota$.
Remember that the Killing form
$$\mathrm{Kil}(X_1, X_2) = \mathrm{trace}(\mathrm{ad}_{X_1} \circ \mathrm{ad}_{X_2}), \qquad
X_1, X_2 \in \mathfrak{g},$$
is non-degenerate when $\mathfrak{g}$ is semisimple.
If $G$ is compact (or more generally, if $G/\mathcal{Z}(G)$ is compact,
where $\mathcal{Z}(G)$ denotes the center of $G$), then
$\mathfrak{g} = \mathfrak{k}$, and $- \mathrm{Kil}(\cdot,\cdot)$,
is positive definite on $\mathfrak{g}$. So, we  can use it to define a bi-invariant Riemannian metric on $G$.
The restriction of this metric to a distribution on $G$
gives a sub-Riemannian manifold. An example of such manifold, $\SU(2)$ (or $S^3$ considered as the set of unit 
quaternions), was  considered in \cite{CChM, ChangMarkVas2}. The problem of geodesic connectivity on $S^3$ was addressed in \cite{ChangMarkVas2}.

If $G$ is non-compact and $\mathfrak{p} \neq 0$ (i.e., $G/\mathcal{Z}(G)$ is non-compact), then
the Killing form restricted to $\mathfrak{p}$ is positive definite.
We consider a left translation of $\mathfrak{p}$ as our horizontal distribution, and using the metric induced by the Killing form, we obtain a sub-Riemannian manifold  equipped with a bi-invariant metric.
Let $K$ be a subgroup of $G$ with the Lie algebra $\mathfrak{k}$. Then $G$ is diffeomorphic to
$K \times \mathfrak{p}$, by $(k, X) \mapsto k e^X$, and if $G$ has a finite center, then $K$ is a maximal compact subgroup. In addition, the above mentioned distribution is horizontal
with respect to the quotient map
$$G \rightarrow G/K.$$
It follows that all normal geodesics are liftings of geodesics from $G/K$, hence they are of the form
$$t \mapsto g_0 e^{t X} e^{- \pr_{\mathfrak{k}} t X} ,\qquad X \in \mathfrak{g},$$
where $g_0 \in G$ is the initial point and $\pr_{\mathfrak{k}}:\mathfrak{g} \to \mathfrak{k}$ is the projection.

Here we consider an example with $G = \SU(1,1)$ and $K = \mathrm{U}(1)$.
Although non-holonomic geometry on $\SU(1,1)$ (or the isometric case of  $\SL(2)$) was first considered earlier in, e.g., \cite{CalinChang, VershikGersh}, we will obtain new results and a more complete description of geodesics both in sub-Riemannian and sub-Lorentzian settings.
Much of all meaningful results come from the analysis of the universal cover of $\SU(1,1)$, which we denote by $\wSU(1,1)$, and which is of its own interest as a new representative of a non-nilpotent Lie group over the topological space $\mathbb R^3$.
We remark also that the K\"ahler manifold  $\SU(1,1)/\U(1)$, is of particular importance in quantum field theories describing  black holes in two-dimensional spacetime by means of an $\SU(1,1)/\U(1)$-gauged Wess-Zumino-Novikov-Witten model, see e.g.,~\cite{Car, Witten}.

Let $\mathfrak{su}(1,1)$ be the Lie algebra of $\SU(1,1)$. When considered as a bilinear form on the entire
$\mathfrak{su}(1,1)$, the Killing form is an index 1 pseudo-metric. Furthermore, the induced Lorentzian metric
 on the Lie group  makes $\SU(1,1)$ isometric to what is called the anti-de Sitter space $\AdS_3$ in General Relativity.
This makes it tempting to  study sub-Lorentzian geometry on $\SU(1,1)$. 
Apart from the fact that the Hamiltonian approach was proposed in \cite{ChangMarkVas}  to study
sub-Lorentzian structures on $\SU(1,1)$, the authors are not aware of  other concrete examples so far, 
where sub-Lorentzian geometry is studied on any other manifold different from the Heisenberg group \cite{Groch1, Groch2} or its extension
to the  $\mathbb H$-type Carnot groups \cite{Korolko}.

The notion of distance in sub-Lorentzian geometry, as well as in Lorentzian geometry,  is given by the supremum of length
over timelike curves. Since timelike loops may appear in $\SU(1,1)$, the distance function behaves badly
(more specifically the distance from a point to itself is $\infty$). Therefore, 
 sub-Lorentzian geometry on $\wSU(1,1)$ is more interesting and meaningful than on $\SU(1,1)$. We are also interested in comparison  of
sub-Lorentzian and Lorentzian geometries on $\wSU(1,1)$, in order to understand somewhat more, how the geometry changes when the Lorentzian metric is being restricted to a distribution. The strong interplay
between $\wSU(1,1)$ and $\SU(1,1)$ will, however, be practical for all explicit calculations.

The structure of the paper is as follows. Section 2 is devoted to general definitions and relations between sub-Riemannian and sub-Lorentzian geometry, and optimal control. Section 3 describes the Lie
groups $\SU(1,1)$ and $\wSU(1,1)$. Section 4 contains results concerning sub-Riemannnian geometry.
We describe the number of geodesics connecting two points and give explicit formulas for
the distance functions. The cut  and conjugate loci on both Lie groups are given.
We discuss the connection between the conjugate locus and the behavior of the geodesics.
In section 5, we completely describe the two-point connectivity problem by sub-Lorentzian geodesics, and compare it with
the  Lorentzian case.  The Lorentzian and sub-Lorentzian future for $\wSU(1,1)$ are compared. It turns out that the reachable sets for Lorentzian and sub-Lorentzian normal geodesics intersect but are not included one to the other.  
Section 6 contains the proofs of main results.

The authors express special thanks to Mauricio Godoy and Irina Markina for many interesting discussions during seminars at the University of Bergen.

\section{Sub Riemannian and Sub-Lorentzian geometry}

\subsection{Sub-Riemannian manifolds}
A sub-Riemannian manifold is an $n$-dimensional manifold $M$, with a fiber metric $\rho$ on an $m$-dimensional smooth distribution
$D \rightarrow M$ ($2 \leq m \leq n$). By distribution, we mean a sub-bundle of the tangent bundle.
Absolutely continuous curves that are almost everywhere tangent to $D$ are called horizontal.
The length of a horizontal curve $\gamma: [0,\tau] \rightarrow M$,
is defined by
$$\ell(\gamma) := \int_0^\tau \rho^{1/2} (\dot{\gamma}(t),\dot{\gamma}(t))
 \, dt.$$
The  {\it Carnot-Carath\'eodory distance} between points $q_0,q_1 \in M$ is defined as
$$d(q_0;q_1) = \inf\limits_{\gamma} \ell(\gamma),$$
where the infimum is taken over all  horizontal curves $\gamma$ satisfying $\gamma(0)=q_0$
and $\gamma(\tau) = q_1$. If there are no such curves connecting $q_0$ and $q_1$, then the
distance is $\infty$.
If  the minimum in the above relation is attained, then the curve is called a length minimizing curve.

Define $D^1 = D$ and iteratively $D^i = D^{i-1} + [D,D^{i-1}]$, $i\geq 2$, where $[\cdot,\cdot]$ denotes the Lie
brackets,  $[X,Y] = XY - YX$ for $X,Y\in T_qM$, $q\in M$. If there exists  a positive integer $k\geq 2$, such that $D^k = TM$, then the distribution is called bracket
generating. It is called regular if $\dim D^i_q$ is independent of the choice of $q$ for all $i$.
We say that $D$ is step $k$ regular, if it is regular and $k$ is the smallest number for which $D^k = TM$. The
Chow-Rashevski{\u\i} theorem \cite{Rashevsky, Chow} states that a bracket generating distribution $D$ guarantees  that any two points of $M$
may be connected by a horizontal curve. In addition, we have the following
generalizations of the corresponding properties from the Riemannian case.

\begin{theorem}[Hopf-Rinow theorem for sub-Riemannian manifolds \cite{BelsR}]
Suppose that $D$ satisfies the bracket generating condition. Then,
\begin{itemize}
\item[i)] any $q_0 \in M$ has a neighborhood $U$, such that there exists a minimizing curve joining  the points $q_1$ and $q_0$ for every $q_1 \in U$;
\item[ii)] if $M$ is complete regarding to $d$, then any two points can be joined by a minimizing curve.
\end{itemize}
\end{theorem}
Observe, that the length minimizing curve may be singular for arbitrarily close points (see \cite{LiuSussman}).

A curve is called geodesic if it is locally a length minimizer.  By a normal geodesic in the sub-Riemannian case we mean  an integral curve of the Hamiltonian system generated by a Hamiltonian function
$H = \sum_{j = 1}^m h_j^2$ in some neighborhood of a point $q$  with respect to
any local orthonormal basis $\{X_1,\dots,X_m\}$ in this neighborhood. Here, and in rest of the paper, if $X$ is a vector field, then $h_X: T^*M\rightarrow \mathbb{R}$ denotes the Hamiltonian function
$$h_X: \lambda \mapsto \lambda(X(q)), \qquad \lambda \in T_q^*M,$$
with respect to the vector field $X$. 
If we have a basis $\{X_1, \dots , X_n \}$, then we simplify notations by writing $h_j$ instead of $h_{X_j}$. 
One of the principal differences between sub-Riemannian  and 
Riemannian geometries is that the function $q \mapsto d(q_0;q)$ is not differentiable
in any neighborhood of $q_0$ when $m<n$.

\subsection{Sub-Lorentzian manifolds}
A sub-Lorentzian manifold is defined  similarly to a sub-Riemannian manifold,
but with $\rho$ now being an index 1 pseudo-metric on $D$.
We will say that a vector $v \in D$ is
\begin{itemize}
 \item timelike if $\rho(v,v) < 0$,
 \item lightlike or null if $\rho(v,v) = 0$,
 \item spacelike if $\rho(v,v) > 0$,
 \item causal or nonspacelike if $\rho(v,v) \leq 0$.
\end{itemize}
A chosen vector field $T$ in $D$, is said to be the time-orientation on $M$, \linebreak if $\rho(T(q),T(q)) < 0$
for any $q \in M$.
A causal vector $v$ is called future directed if $\rho(T(q),v) < 0$, and past directed if  $\rho(T(q),v) > 0$.
A  horizontal curve $\gamma:[0,\tau]\ \rightarrow \mathbb{R}$ is called timelike, null,
spacelike, causal, future directed or past directed,
respectively, if $\dot{\gamma}(t)$ is such a vector for almost every $t \in [0,\tau]$.
We define the timelike future  $\mathcal{I}^+(q_0,\rho)$ of $q_0$ with respect to $\rho$ as the set of all points $q_1 \in M$, such that there is a horizontal,
timelike future directed curve $\gamma$, with $\gamma(0) = q_0$ and $\gamma(\tau) = q_1$.
The causal future, $\mathcal{J}^+(q_0,\rho)$, is defined similarly,  with timelikeness interchanged with causality.
Analogous definitions are valid for the timelike or causal past, which we denote by $\mathcal{I}^-(q_0,\rho)$
and $\mathcal{J}^-(q_0,\rho)$ respectively.
We define the length of a horizontal causal curve by
$\ell(\gamma) = \int_0^\tau |\rho(\dot{\gamma}(t),\dot{\gamma}(t))|^{1/2} dt$.
The sub-Lorentzian distance is  defined by
$$d(q_0;q_1) = \left\{ \begin{array}{ll}
\sup\limits_{\gamma} \ell(\gamma), & \text{if \,} q_1 \in \mathcal{J}^+(q_0,\rho), \\
0, & \text{otherwise}. \end{array} \right.$$
The supremum is taken over all horizontal future directed causal curves from $q_0$ to $q_1$.
Similarly to the Lorentzian distance, the sub-Lorentzian distance satisfies the reverse triangle inequality,
and may not be very well behaving. For instance, if there is a timelike loop trough a point $q \in M$,
then $d(q; q) = \infty$.

A curve $\gamma:[0,\tau] \rightarrow \mathbb{R}$ is called a  length maximizer, if
$\ell(\gamma) = d(\gamma(0);\gamma(\tau))$.
Similarly, a curve $\gamma$ is called a relative maximizer with respect to an open set $W$,
if $\gamma([0,\tau]) \subseteq W$ and $\ell(\gamma) = \sup_{\tilde \gamma}\ell(\tilde \gamma)$,
where the supremum is taken over all horizontal future directed causal curves contained in $W$,
connecting $\gamma(0)$ and $\gamma(\tau)$.
By using the maximum principle, for $D$ bracket generating, we know that all relative maximizers
are either normal geodesics or strictly abnormal maximizers \cite{Groch2}, and that
the relative maximizers always exist locally.
By normal sub-Lorentzian geodesics $\gamma$ we mean curves, such that for any local orthonormal basis $\{X_1,\dots,X_m\}$ of $D$,
with $X_1$ as the time-orientation, $\gamma$ is an integral curve of  the Hamiltonian system generated by a Hamiltonian function $-h_1^2 + \sum_{j=2}^m h_j^2$.

The question of whether length maximizers exist between two points, is a much more complicated
than the question of the existence of length minimizers in Riemannian or sub-Riemannian geometries.
The most common sufficient condition for the global existence of maximizing curves on a
Lorentzian or sub-Lorentzian
manifold $M$ is a relatively strict requirement that $M$ should be globally hyperbolic, i.e.,
strongly causal (every point has an arbitrarily  small neighborhood, such that causal curves that leave the neighborhood never return back), and that $\mathcal{J}^+(q_1,\rho) \cap \mathcal J^{-}(q_2,\rho)$ is compact for
any $q_1,q_2 \in M$.

\subsection{Minimizing and maximizing curves seen from the viewpoint of optimal control}
Determination of  curves whose length is equal to the distance, either in the
sub-Riemannian or sub-Lorentzian setting, can be formulated as a solution to an optimal control problem.
Let $M$ be an $n$-dimensional manifold, let a submersion $\pi: \,V \rightarrow M$ be a fiber bundle with the fiber $U=\pi^{-1}(q)$, $q\in M$,
and let $f:V \rightarrow TM$ be a fiber preserving map.

An {\it  admissible pair} $v=(u,\gamma_u)$ is a bounded measurable mapping
$v:[0,\tau] \rightarrow V, 0 < \tau < \infty$ such that $\gamma_u = \pi \circ v$ is a Lipschitzian curve
on $M$ and $\dot{\gamma}_u = f(u, \gamma_u)$. The component $u$ is called {\it control} in the literature.
For some fixed value of $\tau$, we denote the space of all admissible pairs by $\mathcal{V}_{\tau}$. The space
$\mathcal{V}_{\tau}$ is a smooth Banach submanifold of  $L^{\infty}([0,\tau],V)$, which is  modeled on $L^\infty([0,\tau],\mathbb{R}^{\dim V})$.
An admissible pair is said to connect $q_0,q_1 \in M$, if $\gamma_u(0) = q_0$, $\gamma_u(\tau) = q_1$.
Such a system is called {\it controllable}, if any two points in $M$ are connected by at least one admissible pair.
Let $J:\mathcal{V}_{\tau} \rightarrow \mathbb{R}$ be a functional, given by
$$J(v) = \int_0^{\tau} \mathcal{K}(v) dt \, ,$$
where $\mathcal{K}$ is some differentiable function. Suppose we are given two points $q_0, q_1 \in M$,
such that there is at least one admissible pair  connecting them.
An {\it optimal control problem} with respect to the functional $J$,
is a problem of finding an element $v^*=(u^*,\gamma_{u^*}) \in \mathcal V_\tau$ connecting $q_0,q_1$, 
such that for any $v \in  \mathcal{V}_\tau$ connecting the same two points,
we have $J(v^*) \leq J(v)$. Analogously, we may look for the maximum of $J$. We call $u^*$ the optimal control, and $\gamma_{u^*}$ the optimal trajectory.
We may also consider {\it free-time} optimal control problems where $\tau$ is allowed to vary.

We will write $\vec{H}$ to denote the Hamiltonian vector field associated to a Hamiltonian function $H$.
For a pseudo-Hamiltonian function $\mathscr{H}$, the vector field $\vec{\mathscr{H}}$ is defined so that for a fixed $u \in U$,
$\vec{\mathscr{H}}(u, \cdot )$ is the Hamiltonian vector field associated with $\mathscr{H}(u, \cdot)$.
The main tool to solve optimal control problems is the following first order condition known as the  {\it Pontryagin Maximum Principle} (PMP).
\begin{theorem}[PMP for Optimal Control Problem with fixed time $\tau$]
For a given  value of $\tau$, let $v^*=(u^*,\gamma_{u^*}) \in \mathcal{V}_\tau$ be an optimal pair for  the above problem, i.e.,
$J(v^*) = \min\limits_{v\in \mathcal{V}_\tau} \{ J(v) \big| \, \gamma_u(0) = 0, \gamma_u(\tau) = q_1 \}$.
For each $\varkappa\in \mathbb R$ a pseudo-Hamiltonian function $\mathscr{H}_{\varkappa}:\, U\times T^*M\to \mathbb R$ is defined by
$$\mathscr{H}_{\varkappa}(u,\lambda) = \lambda(f(u,q)) + \varkappa \mathcal{K}(u,q), $$
where $\lambda \in T_q^*M$.
Then there exists a curve $\xi:[0,\tau] \rightarrow T^*M$, and a number $\varkappa \leq 0$, such
that 
$$\dot{\xi}(t) =\vec{\mathscr{H}}_{\varkappa,}(u^*,\xi(t)), \qquad
\mathscr{H}_{\varkappa}(u^*,\xi(t)) = \max_{u \in U} \mathscr{H}_\varkappa(u,\xi(t)) 
\quad \text{a.e. in} \, \, t \in [0,\tau],$$
and
$$\pr_M \xi(t) = \gamma_{u^*}(t).$$
The pseudo-Hamiltonian function is constant along the optimal trajectory $\xi(t)$.
Moreover, if $\varkappa = 0$, then for almost every $t\in [0,\tau]$, $\xi(t)$ can not be in the zero section of $T^*M$.
\end{theorem}

For the problem of the maximum of  $J$, the above theorem has the same formulation, changing only the requirement  $\varkappa \leq 0$ to $\varkappa \geq 0$.
If we consider a free-time problem, then we also require $\mathscr{H}_{\varkappa}(u^*, \xi(t)) \equiv 0$ almost everywhere.
If $\varkappa \neq 0$, then the solution is called normal (in this case we may just choose $\varkappa = \pm 1$).
If $\varkappa = 0$, the solution is called abnormal.

\begin{remark}[\cite{AS}]
If $H_\varkappa(\lambda) = \max_{u \in U} \mathscr{H}_\varkappa(u,\lambda)$ is defined and is $C^2$ on
$T^*M \setminus s_0(M)$, where $s_0: M \rightarrow T^*M$ is the zero-section, then
$\dot{\xi}(t) = \vec{H}_\varkappa(\xi(t))$.
\end{remark}

Let us turn to  specific sub-Riemannian and sub-Lorentzian settings.
We only need the formulation for the case when $M = G$ is a Lie group.
This case is somewhat simpler since the tangent bundle of Lie groups is trivial, so
we can always find a global basis to span a distribution.
Let $D = \mathrm{span}\{X_1, \dots, X_m\}$ be a distribution, and $\rho$ be either a sub-Riemannian or sub-Lorentzian
metric on $D$. We may assume that the collection $\{X_j\}_{j=1}^m$ forms an orthonormal basis with respect to the metric $\rho$.
Let $V = U \times G$, where $U \subset \mathbb{R}^m$.
We choose $f(u_1,\dots,u_m, g) = u_1 X_1(g) + \cdots + u_m X_m(g)$,  $g\in G$, to
ensure that $g(t)$ is a horizontal curve.

The determination of the sub-Riemannian distance, comes down to finding
an optimal pair $(u^*, g_{u^*}):[0 , \tau] \rightarrow U \times G$
that minimizes the functional
$$J(u,g_u) = \frac{1}{2} \int_0^{\tau}  \sum_{j=1}^m u_j^2 \, dt, \, 
\quad u \in U = \{u \in \mathbb{R}^m | \, \|u\| \leq 1 \}.$$
The corresponding pseudo-Hamiltonian is given by
$$\mathscr{H}_\varkappa (u,\lambda) = -\frac{\varkappa}{2} \sum_{j = 1}^m u_j^2
+ \sum_{j= 1}^m u_j h_j(\lambda), \quad \lambda \in T^*G, \quad u \in U, \quad  \varkappa = 0,1. $$
As a consequence of Fillipov's theorem \cite{AS}, there always exists
a length minimizing horizontal curve between two points in this setting. 

Similarly, finding the sub-Lorentzian distance with $X_1$ as a time-orientation,
 means finding a pair $(u^*,g_{u^*})$ that maximizes
$$J(u,g_u) = \int_0^{\tau} \left(u_1^2 - \sum_{j=2}^m u_j^2\right)^{\frac{1}{2}} \, dt, \quad 
\qquad U = \left\{u \in \mathbb{R}^m \big| \, u_1 = \left(1+ \sum_{j =2}^m u_j^2\right)^{\frac{1}{2}}
\right\},$$
and the pseudo-Hamiltonian related to this problem becomes
$$\mathscr{H}_\varkappa(u,\lambda) = \varkappa\left(u_1^2 - \sum_{j=2}^m u_j^2\right)^{\frac{1}{2}}
 + \sum_{j= 1}^m u_j h_j(\lambda).$$
The latter case is  more complicated because we optimize a non-convex functional.

In our example of $\SU(1,1)$, we only consider left-invariant distributions, that
lead to left-invariant Hamiltonians. Let us, therefore, include some sketch of the theory of 
Hamiltonian systems in this case.
Generally, assume that $H$ is a left-invariant Hamiltonian function on a Lie group
(i.e., $H(d^*L_g \lambda) = H(\lambda)$). From the isomorphisms of bundles, we have
$$G \times \mathfrak{g} \rightarrow TG, \qquad G \times \mathfrak{g}^* \rightarrow T^* G,$$
$$(g,X) \mapsto dL_g X, \qquad (g,p) \mapsto d^*L_{g^{-1}} p.$$
The traditional Hamiltonian equation $\dot{\lambda}(t) = \vec{H}(\lambda(t))$, which locally has the
form
$$\dot{g} = \frac{\partial H}{\partial \xi}, \qquad \dot{\xi} = - \frac{\partial H}{\partial g}, \qquad
\lambda = (g,\xi) \in T^*M,$$
may be written as
$$\dot{g} = dL_g \left(\frac{\partial H}{\partial p}(g,p) \right),$$
$$\dot{p} = d^*L_g \left( \frac{\partial H}{\partial g}(g,p) \right) +
\mathrm{ad}\left(\frac{\partial H}{\partial p}(g,p) \right)^* (p).$$

If $H$ is left-invariant, and hence is independent of $g$, then the first term in the
equation for $\dot{p}$ vanishes. In particular, if $X$ is a left-invariant vector field, then
$h_X(g,p) = d^*L_{g^{-1}}p(X(g)) = p(dL_{g^{-1}} dL_g X(1)) = p(X(1))$, so
$h_j(p) = p(X_j) =: p_j$.
Furthermore,
\begin{eqnarray*}\dot{p}_j &=& \dot{p}(X_j) = \mathrm{ad}\left(\frac{dH}{dp}(p)\right)^*(p)(X_j)
= p\left(\left[ \frac{dH}{dp}(p) , X_j \right] \right) \\ &=& p\left(\left[ \frac{dH}{dp}(p) , \frac{dh_j}{dp}(p) \right] \right)
= \{h_j , H \}(p).
\end{eqnarray*}
Here $\{ \cdot , \cdot \}$  denotes the Poisson bracket.
For more details, see \cite{Jur1} (observe that the difference in sign in our formulation and in \cite{Jur1} happens
because of  different definitions of $[ \cdot , \cdot ]$.)
Notice that similar considerations can be done by interchanging the left and the right actions.

\section{Structure of the Lie groups $\SU(1,1)$ and $\wSU(1,1)$}

\subsection{Some notation} \label{notationsec}
We will start by defining the following functions, which will simplify our notations later in dealing with the exponential map.
For any $(\alpha,t) \in \real \times \real_+$, define
\begin{equation} \label{defscrCscrS} \nonumber \scrC_\alpha(t) := \left\{ \begin{array}{ll} \cosh(\sqrt{\alpha} \cdot t) & \text{if } \alpha \geq 0 \\
\cos(\sqrt{-\alpha} \cdot t) & \text{if } \alpha < 0 \\ \end{array} \right. ,  \qquad
\scrS_\alpha(t) := \left\{ \begin{array}{ll} \frac{\displaystyle \sinh(\sqrt{\alpha} \cdot t)}{\displaystyle\sqrt{\alpha}}& \text{if } \alpha > 0 \\
t & \text{if } \alpha = 0 \\
\frac{\displaystyle\sin(\sqrt{-\alpha} \cdot t)}{\displaystyle\sqrt{-\alpha}} & \text{if } \alpha < 0 \\ \end{array} \right. . \end{equation}

These functions are related by the equation $\scrC_\alpha^2(t) - \alpha \scrS_\alpha^2(t) = 1.$ One can check that the partial derivatives of these functions satisfy the following relations.
\begin{equation} \label{derivativesCS} \nonumber
\begin{array}{ccc} \frac{\partial}{\partial t} \scrC_\alpha(t) = \alpha \scrS_\alpha(t), & \qquad &
\frac{\partial}{\partial t} \scrS_\alpha(t) =  \scrC_\alpha(t), \\ \\
\frac{\partial}{\partial \alpha} \scrC_\alpha(t) = \frac{t}{2} \scrS_\alpha(t), &
\qquad & \frac{\partial}{\partial \alpha} \scrS_\alpha(t) = \frac{t}{2 \alpha}\left(\scrC_\alpha(t) - \scrS_\alpha(t)\right).\end{array}\end{equation}
Observe that $\frac{\partial}{\partial \alpha} \scrS_\alpha(t)|_{\alpha = 0}$ does not exist, because
$$ \lim_{\alpha \rightarrow 0 +}\tfrac{\partial}{\partial \alpha} \scrS_\alpha(t) = \tfrac{t}{6},\quad \mbox{while \ }
\lim_{\alpha \rightarrow 0 -}\tfrac{\partial}{\partial \alpha} \scrS_\alpha(t) = -\tfrac{t}{6}.$$
When $t = 1$, we will write $\scrC_\alpha := \scrC_\alpha(1)$ and $\scrS_\alpha := \scrS_\alpha(1)$.

\begin{remark}\label{chi}
$\frac{\partial}{\partial \alpha} \scrC_\alpha = 0$ implies that $\alpha = - \pi^2 k^2$, while $\frac{\partial}{\partial \alpha} \scrS_\alpha = 0$ implies $\alpha = - \chi_k^2$, where $\chi_k$, for each $k \in \mathbb{N}$, is defined as a unique number satisfying
\begin{equation} \label{defchik} \chi_k = \tan \chi_k, \qquad
\pi k \leq \chi_k \leq \pi k + \frac{\pi}{2}.\end{equation}
The numbers $\chi_k$ will become important later.
\end{remark}

The last function we need to define is $\scrT: \real^3 \to \real$.
We do this by  associating any $a = (a_1,a_2,a_3) \in \real^3$ to a real number $\check a$ defined by
\begin{equation} \label{defchecka} \nonumber \check a := a_1^2 + a_2^2 - a_3^2.\end{equation}
If $a$ satisfies $\check a \geq 0$, then we define $\scrT_a$ to be
\begin{equation} \label{defscrT1} \nonumber \scrT_a = \left\{ \begin{array}{ll} 
\tan^{-1} \left(\frac{a_3}{\sqrt{\check a}}\tanh \sqrt{\check a}\right) & \text{if } \check a > 0 \\\
\tan^{-1} a_3 & \text{if } \check a = 0\end{array}\right. . \end{equation}
If $\check a <0$, then $\scrT_a$ is given by
\begin{equation} \label{defscrT2} \nonumber \scrT_a = \left\{ \begin{array}{ll}
\tan^{-1}\left(\frac{a_3}{\sqrt{-\check a}} \tan \sqrt{-\check a}) \right)
+ \mathrm{sgn}(a_3) \pi \left\lceil \frac{\sqrt{-\check a}}{\pi} - \frac{1}{2} \right\rceil & \text{if } \sqrt{-\check a}
\neq \pi/2 \text{ mod } \pi \\
\mathrm{sgn}(a_3) \sqrt{- \check a} & \text{if } \sqrt{-\check a} = \pi/2 \text{ mod } \pi 
\end{array} \right. .\end{equation}
Here, $\lceil x \rceil$ denotes the ceiling of $x$ (i.e.,
$\lceil x \rceil = \min \{j \in \mathbb{Z} \, | \,  j \geq x \}$) and
$$\mathrm{sgn} \, x = \left\{\begin{array}{ll} 1 & \text{if } x> 0 \\ 0 & \text{if } x = 0 \\ -1 & \text{if } x<0 \end{array} \right. .$$
$\scrT$ is connected to the other two functions by the formula
$$e^{i\scrT_a} = \frac{\scrC_{\check a} + i a_3 \scrS_{\check a}}{|\scrC_{\check a} + i a_3 \scrS_{\check a}|}.$$

\subsection{Lie group $\SU(1,1)$}
We consider the Lie group $\SU(1,1)$ of unitary complex $(2\times 2)$ matrices of the form
$$g = \left( \begin{array}{cc} z_1 & z_2 \\ \bar{z}_2 & \bar{z}_1 \end{array} \right),$$ with the operation of
usual matrix multiplication. Simplifying notation, we will sometimes view $g \in \SU(1,1)$ as an element
from $\mathbb{C}^2$, writing $g =(z_1, z_2)$.
In this notation, the group operation is
$$(z_1, z_2) (\hat{z}_1, \hat{z}_2) = \left(z_1 \hat{z}_1 + z_2 \overline{\hat{z}_2},
z_1 \hat{z}_2 + z_2 \overline{\hat{z}_1}\right),$$
the identity is $1 = (1,0)$,  and the inverse element becomes $(z_1,z_2)^{-1} = (\bar{z}_1, -z_2)$.
$\SU(1,1)$ is isomorphic to the Lie group $\SL(2)$ of real $2 \times 2$ matrices of determinant 1, by
the mapping
$$\begin{array}{ccc}\SU(1,1) & \cong & \SL(2), \\
\left( \begin{array}{cc} z_1 & z_2 \\ \bar{z}_2 & \bar{z}_1 \end{array} \right) & \mapsto
& \left( \begin{array}{ccc} \re z_1 + \im z_2 & \, & - \im z_1 - \re z_2 \\ \im z_1 - \re z_2 & \, & \re z_1 - \im z_2
\end{array} \right). \end{array}$$ 

We will often use the isomorphism from the tangent bundle $T\SU(1,1)$ to a sub-bundle of the complex
tangent bundle $T_{\mathbb{C}} \SU(1,1)$, given at any fixed point  by
$$u_1 \partial_{x_1} + u_2 \partial_{y_2} + u_3 \partial_{y_3} + u_4 \partial_{x_4} \mapsto
(u_1 + i u_2) \partial_{z_1} + (u_1 - i u_2) \partial_{\bar{z}_1} + 
(u_3 + iu_4) \partial_{z_2} + (u_3 - iu_4) \partial_{\bar{z}_2},$$
where $z_j = x_j + i y_j$, $\partial_{z_j} = \frac{1}{2} (\partial_{x_j} - i\partial_{y_j})$ and
$\partial_{\bar{z}_j} = \frac{1}{2} (\partial_{x_j} + i\partial_{y_j})$.
Let $L_g$ denote the left actions by $g$. The tangent map of the left action is given by
$$dL_g: b_1 \partial_{z_1} + \bar{b}_1 \partial_{\bar{z}_1}
+ b_2 \partial_{z_2} + \bar{b}_2 \partial_{\bar{z}_1}\mapsto$$
$$\mapsto
(b_1 z_1 + \bar{b}_2 z_2) \partial_{z_1} + (\bar{b}_1 \bar{z}_1 + b_2 \bar{z}_2) \partial_{\bar{z}_1}
+ (b_2 z_1 + \bar{b}_1 z_2) \partial_{z_2} + (\bar{b}_2 \bar{z}_1 + b_1 \bar{z}_2) \partial_{\bar{z}_2}.$$ 
The tangent space at the identity is spanned by the vectors
$$X(1) = - \partial_{z_2} - \partial_{\bar{z}_2}, \quad Y(1) = i\partial_{z_2} - i\partial_{\bar{z}_2},
\quad Z(1) = -i \partial_{z_1} + i \partial_{\bar{z}_1}.$$
Viewed as elements of the Lie algebra $\mathfrak{su}(1,1)$ of $\SU(1,1)$, they have the following form
\begin{equation} \label{Basis}
X = \left(\begin{array}{cc} 0 & -1 \\ -1 & 0 \end{array} \right), \quad
Y = \left(\begin{array}{cc} 0 & i \\ -i & 0 \end{array} \right), \quad
Z = \left(\begin{array}{cc} -i & 0 \\ 0 & i \end{array} \right).
\end{equation}
So we can obtain the left-invariant basis for the tangent bundle as
\[\begin{array}{ll} 
X(g) & = - z_2 \partial_{z_1} - \bar{z}_2 \partial_{\bar{z}_1} - z_1 \partial_{z_2}- \bar{z}_1 \partial_{\bar{z}_2}, \\
Y(g) & = - i z_2 \partial_{z_1} + i\bar{z}_2 \partial_{\bar{z}_1} +  iz_1 \partial_{z_2} - i \bar{z}_1 \partial_{\bar{z}_2}, \\
Z(g) & = - i z_1 \partial_{z_1} + i\bar{z}_1 \partial_{\bar{z}_1} + i z_2 \partial_{z_2} - i \bar{z}_2 \partial_{\bar{z}_2}.
\end{array}\]
We choose a corresponding dual basis for the cotangent bundle, given by
\begin{align*} X^*(g) &= \frac{1}{2}(\bar{z}_2 dz_1 + z_2 d\bar{z}_1 - \bar{z}_1 dz_2-  z_1 d\bar{z}_2), \\
Y^*(g) &= \frac{1}{2}(-i\bar{z}_2 dz_1 + i z_2 d\bar{z}_1 - i\bar{z}_1 dz_2 +  i z_1 d\bar{z}_2), \\
Z^*(g) &= \frac{1}{2}(i\bar{z}_1 dz_1 - iz_1 d\bar{z}_1 + i\bar{z}_2 dz_2-  i z_2 d\bar{z}_2). \end{align*}
The bracket relations yield
$$[X,Y] = - 2 Z, \qquad [X,Z] = - 2 Y, \qquad  [Y,Z] = 2 X.$$
It follows that any distribution spanned by two of three vector fields is bracket generating.

The exponential map for this Lie group is
\begin{align} \label{SU11exp} e^{a_1X + a_2 Y + a_3 Z} &= \left( \begin{array}{ccc}
\scrC_{\check a} - i a_3 \scrS_{\check a}&\, & - (a_1 - ia_2) \scrS_{\check a} \\ 
- (a_1 + ia_2) \scrS_{\check a} & \, & \scrC_{\check a} + i a_3 \scrS_{\check a} 
\end{array} \right) \\ \nonumber
& = \scrC_{\check a} \cdot 1 + \scrS_{\check a} (a_1X + a_2Y + a_3Z).\end{align}

We define the metric $\rho$ as a left-invariant metric, whose restriction to the $\mathfrak{su}(1,1)$
has the following form
$$(X_1,X_2) \mapsto  \frac{1}{8} \mathrm{Kil}(X_1,X_2)$$
Remark that in $\mathfrak{su}(1,1)$, the Killing form is equal to
$$\mathrm{Kil}(X_1,X_2) = 4 \, \mathrm{trace}(X_1X_2).$$
In the basis of $\partial_{x_1}, \partial_{y_1}, \partial_{x_2}, \partial_{y_2}$, the metric tensor of $\rho$
has the form
$$(\rho_{ij}) = \left(
\begin{array}{cccc} -1 & 0 & 0 & 0 \\ 0 & -1 & 0 & 0 \\ 0 & 0 & 1 & 0 \\ 0 & 0 & 0 & 1 \end{array} \right).$$
The metric $\rho$ is then  Lorentzian, and  in fact,  bi-invariant.
Hence, as a Lorentzian manifold, $\SU(1,1)$ may be considered as a subset of
$\mathbb{R}^{2,2}$, which is $\mathbb{R}^4$ with an index 2 pseudo-metric.
This subset 
$$\AdS_3 = \left\{(x_1, y_1, x_2, y_2) \in \mathbb{R}^{2,2}| \, x_1^2 + y_1^2 - x_2^2 - y_2^2 = 1 \right\}$$
is called the 3-dimensional anti-de Sitter space.

The restriction of $\rho$ to the distribution $D = \spn\{X,Y\}$, makes it positive definite, and makes
$(\SU(1,1),D, \rho|_D)$ a sub-Riemannian manifold. Similarly, the restriction $\rho$ to
$E = \spn\{Y,Z\}$ makes $(\SU(1,1),E,\rho|_E)$ a sub-Lorentzian manifold, if we define the
time-orientation by $Z$.
The latter case however, contains timelike loops. An example is
$$t \mapsto \left(\begin{array}{cc} e^{-it} & 0 \\ 0 & e^{it} \end{array} \right),$$
which is a loop through 1. To avoid this problem we  will also study the universal cover
$\wSU(1,1)$ of $\SU(1,1)$. This will also be helpful in order to study sub-Riemannian geometry
on $\SU(1,1)$. In addition, it is an interesting example on its own. In fact, sometimes  the attribution
anti-de Sitter space is used for the universal cover instead of $\AdS_3$ itself.

\subsection{Lie group structure of the universal cover $\wSU(1,1)$ of $\SU(1,1)$}
Since $\SU(1,1)$ is diffeomorphic to $S^1\times  \mathbb{R}^2$, the universal cover must be diffeomorphic to
$\mathbb{R}^3$.
We represent the covering space $\wSU(1,1)$
as $\mathbb{R}\times \mathbb{C}$, with the covering map $\pi:\wSU(1,1) \rightarrow \SU(1,1),$ defined by
\begin{equation} \label{projwSU11toSU11} \tilde g = (c,w) \mapsto \left( \begin{array}{ccc} \sqrt{1 + |w|^2} e^{ic} &\, & w \\ \bar{w} &\, & \sqrt{1 + |w|^2} e^{-ic}
\end{array} \right),\end{equation}
where $c \in \mathbb{R}$ and $w \in \mathbb{C}$.
We define the product on $\wSU(1,1)$, to be the unique product for which
 $\tilde 1 = (0 , 0)$ is the identity, and which makes $\pi$ a Lie group homomorphism.
It is obvious that the Lie algebra of $\wSU(1,1)$ is also $\mathfrak{su}(1,1)$.

\begin{definition}
Let $(c_j, w_j) \in \wSU(1,1), \, \, j = 1,2$. We define the operation on $\wSU(1,1)$ as $(c_1,w_1)(c_2,w_2) = (c,w)$, where
$$c =c_1 + c_2 + \tan^{-1}\left(
\frac{\im(w_1\bar{w}_2 e^{-i(c_1 +c_2)})}
{\sqrt{(1 +|w_1|^2)(1+|w_2|^2)}  + \re(w_1\bar{w}_2 e^{-i(c_1 +c_2)})}\right),$$
$$w = w_2 \sqrt{1+|w_1|^2} e^{ic_1} + w_1\sqrt{1+|w_2|^2} e^{-ic_2}.$$
\end{definition}

\begin{proposition} \label{groupop}
The above definition provides a group structure to $\wSU(1,1)$ with the identity $(0,0)$. With respect to this group structure, $\pi: \wSU(1,1) \to \SU(1,1)$ becomes a Lie group homomorphism.
\end{proposition}

\begin{proof}
The fact that $\pi(c_1,w_1) \pi(c_2,w_2) = \pi(c,w)$, trivially implies the expression for
$w$. The value of $c$ must satisfy the relation
$$e^{ic} = \exp\left(i \Arg\left(\sqrt{(1+|w_1|^2)(1+|w_2|^2)} e^{i(c_1+c_2)} + w_1 \bar{w}_2\right) \right)$$
$$= e^{i(c_1 + c_2)} \exp\left(i \Arg\left(\sqrt{(1+|w_1|^2)(1+|w_2|^2)} + w_1 \bar{w}_2e^{-i(c_1+c_2)} \right) \right).$$
Now since
$$\left|\sqrt{(1+|w_1|^2)(1+|w_2|^2)}\right| > \left|w_1 \bar{w}_2e^{-i(c_1+c_2)}\right|,$$
we know that
$$-\frac{\pi}{2} < \Arg \left(\sqrt{(1+|w_1|^2)(1+|w_2|^2)} + w_1 \bar{w}_2e^{-i(c_1+c_2)}
\right) < \frac{\pi}{2},$$
and the formula for $c$ follows. It is clear that $(0,0)$ is the identity under this product.
Observe that $(c,w)^{-1} = (-c, - w)$.
The associativity of the product remains to be proved.
Due to the associativity of the product on  $\SU(1,1)$, we only need to show that if
$$(C,W) = \Big((c_1,w_1) (c_2,w_2) \Big) (c_3,w_3), \qquad
(\widehat{C},\widehat{W}) = (c_1,w_1)\Big( (c_2,w_2) (c_3,w_3) \Big),$$
then $0 = C - \widehat{C}$. Let $\pi(C,W) = \pi(\widehat{C},\widehat{W}) = (Z_1,Z_2)$, and denote by
$$\theta = \Arg\left(\sqrt{(1+|w_1|^2)(1+|w_2|^2)} + w_1 \bar{w}_2 e^{-i(c_1 +c_2)}\right),$$
$$\widehat{\theta} = \Arg\left(\sqrt{(1+|w_2|^2)(1+|w_3|^2)} + w_2 \bar{w}_3 e^{-i(c_1 +c_2)}\right).$$
Then
$$C = c_1 + c_2 + c_3 + \theta + \Arg(Z_1 e^{-i(c_1+c_2 + c_3 + \theta)}), \quad
\widehat{C} = c_1 + c_2 + c_3 + \widehat{\theta} + \Arg(Z_1 e^{-i(c_1+c_2 + c_3 + \hat{\theta})}).$$
Similarly to the reasoning above, $$ \Arg(Z_1 e^{-i(c_1+c_2 + c_3 + \theta)})\,\,\,\mbox{and}\,\,\,
\Arg(Z_1 e^{-i(c_1+c_2 + c_3 + \widehat{\theta})}) \in (  -\frac{\pi}{2},\,\frac{\pi}{2}),$$ so since
$|\theta - \widehat{\theta}| < \pi$,
we know that $\Arg(Z_1 e^{-i(c_1+c_2 + c_3 + \theta)})- \Arg(Z_1 e^{-i(c_1+c_2 + c_3 + \hat{\theta})})
= \widehat{\theta} - \theta$, and it follows that $C = \widehat{C}$.
\end{proof}

The left action is given by ($b_0 \in \mathbb{R}, b \in \mathbb{C}$)
$$dL_{\tilde g} : b_0 \partial_c + b \partial_w + \bar{b} \partial_{\bar{w}}
\mapsto$$
$$(b_0 +  \im(\bar{b} w e^{-ic}))\partial_c + (-ib_0w + b \sqrt{1 + |w|^2} e^{ic}) \partial_w
+ (ib_0 \bar{w} + \bar{b} \sqrt{1 + |w|^2} e^{-ic}) \partial_{\bar{w}}$$
$$= \left(b_0 +  \im\left(\frac{\bar{b} z_2|z_1|}{z_1}\right)\right)\partial_c + (-ib_0z_2 + b z_1) \partial_w
+ (ib_0 \bar{z}_2 + \bar{b} \bar{z}_1) \partial_{\bar{w}}.$$
The vector fields $X,Y,Z$ lifted to the universal cover become
$$\widetilde X(c,w) = - \im(we^{ic}) \partial_c + \sqrt{1 + |w|^2} (e^{ic} \partial_w + e^{-ic}\partial_{\bar{w}}),$$
$$\widetilde Y(c,w) = -\re(we^{ic}) \partial_c -i \sqrt{1 + |w|^2} (e^{ic} \partial_w - e^{-ic}\partial_{\bar{w}}),$$
$$\widetilde Z(c,w) = - \partial_c + iw \partial_w -i \bar{w} \partial_{\bar{w}}.$$
The exponential map in $\wSU(1,1)$ is
\begin{equation} \label{wSU11exp} e^{a_1X + a_2Y + a_3Z} = \big(- \scrT_a, - (a_1 - ia_2) \scrS_{\check a}\big),\end{equation}
which can be found by lifting the exponential map from $\SU(1,1)$.

We lift the metric $\rho$ to a metric on $\wSU(1,1)$. Denote this metric by $\tilde \rho$.
Analogously to $\SU(1,1)$, let us define the distributions
$\widetilde D = \spn\{ \widetilde X, \widetilde Y \}$ and
$\tilde E = \spn\{ \widetilde Y, \widetilde Z \}$. The restriction of $\tilde \rho$ to $\widetilde D$ and $\widetilde E$,
defines  sub-Riemannian and sub-Lorentzian structures on $\wSU(1,1)$ respectively.
In the sub-Lorentzian case, we let $\widetilde Z$ define the time orientation.

\begin{remark}
We can construct a sub-Lorentzian manifold, by considering the distribution
$\spn\{\widetilde X,\widetilde Z\}$, but the geodesics are very similar to the sub-Lorentzian manifold \linebreak
$(\wSU(1,1), \tilde E, \widetilde \rho|_{\tilde E})$ so we omit this choice of distribution.
\end{remark}

\section{Sub-Riemannian geometry on $\SU(1,1)$ and $\wSU(1,1)$}

\subsection{Geodesics,  horizontal space, and  vertical space}
We will now take advantage of the fact that the pseudo-metric induced by the Killing form is bi-invariant.
\begin{theorem} \label{thmon}
Let $G$ be a Lie group with the Lie algebra $\mathfrak{g}$, and with a bi-invariant pseudo-metric $\rho$.
Let $K$ be a subgroup of $G$, with the Lie algebra $\mathfrak{k}$, and let us denote $\mathfrak{p}
= \mathfrak{k}^{\perp}$. Define a left-invariant distribution $D$ by $D_g = dL_g \mathfrak{p}$.
Then all normal geodesics on the non-holonomic manifold $(G, D, \rho|_D)$
are lifting of the normal geodesics on $G/K$ with the induced metric. This means that
all normal geodesics starting at $g_0 \in G$ are of the form
$$t \mapsto g_0 e^{X_0t} \cdot e^{- \pr_\mathfrak{k} X_0 t}, \qquad X_0 \in \mathfrak{g} \,,$$
where $\pr_{\mathfrak{k}}:\mathfrak{g} \rightarrow \mathfrak{k}$ is the projection. 
\end{theorem}
This theorem is a special case of the corresponding result from \cite{Mon}.
We may use Theorem \ref{thmon}  with $K = e^{\spn \{Z\}}$.
The normal geodesics starting at the identity, by equations \eqref{SU11exp} and \eqref{wSU11exp}, admit the form
$$z_1(t) = (\scrC_{\check a}(t) - i a_3 \scrS_{\check a}(t)) e^{ia_3t}, \qquad z_2(t) = - (a_1 - ia_2) \scrS_{\check a}(t) e^{-ia_3t},$$
in $\SU(1,1)$ and
$$c(t) = - \scrT_{t a} + a_3 t, \qquad w(t) = - (a_1 - ia_2) \scrS_{\check a}(t) e^{-ia_3t},$$
in $\wSU(1,1).$ See Figure 1 for examples. Any other normal geodesic starting at some point, is a left translation of a normal one starting at the identity.

Remark that if $a_3 = 0$ (i.e. $\pr_{\mathfrak k} X_0 = 0$), then the sub-Riemannian geodesics become just curves given by the exponential map. We define {\it the horizontal space} as the collection of points that can be reached by such geodesics. In $\wSU(1,1)$ this is the collection of points $\tilde g = (0,w)$.
In $\SU(1,1)$ these points are on the form $g = (z_1, z_2), z_1 >1$. We will also use the term {\it vertical space} or {\it vertical line} for the points in $e^{\spn \{Z\}}$ in both $\SU(1,1)$ and $\wSU(1,1)$.

\begin{figure}[ht!]
\centering
\scalebox{0.50}{\includegraphics{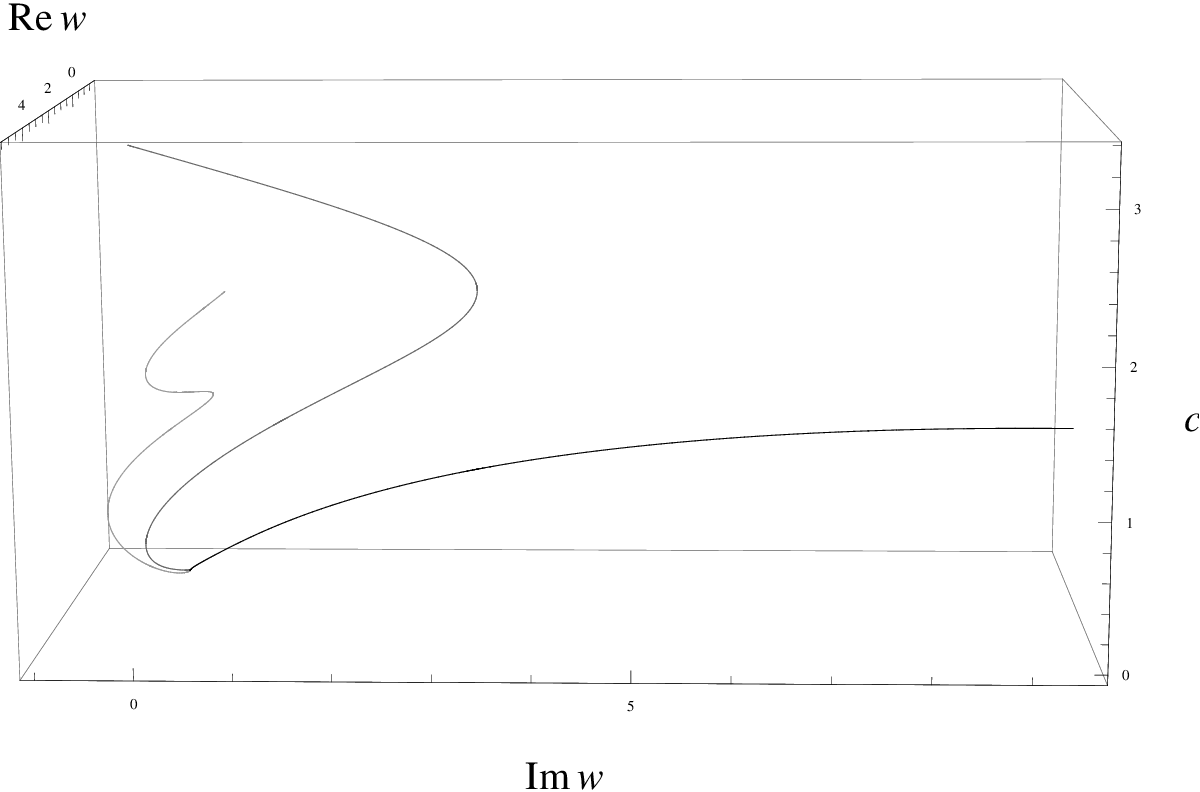}}
\caption[]{The above figure shows three geodesics in $\wSU(1,1)$, all of length $\frac{3 \pi}{2}$.
The figure shows geodesics starting at $\tilde 1$, which, when parametrized by arc length,
have initial values of $\check a$ given (from left to right) by $-1, 0$ and $\frac{3}{4}$.}
\end{figure}

\subsection{Length and number of geodesics}
Recall the definition of $\chi_k$ in \eqref{defchik}.
Then, up to reparameterization, we have the following results 
regarding the number of geodesics connecting $\tilde 1$ and $\tilde g \in \wSU(1,1)$.

\begin{proposition} \label{sRnumbergeo}
\begin{itemize}
\item[(a)] If $\tilde g = (c,0)$ for some $c$, then the number of
geodesics is uncountable (countably many geometrically different).
\item[(b)] If $\tilde g = (0,w)$ for some $w$, then there exists a unique geodesic connecting
$\tilde g$ and $\tilde 1$. This geodesic is contained in the plane $c= 0$.
\item[(c)] If $\tilde g$ is any other point, then we can obtain the number of geodesics in the following way. Let $k$ be the largest positive integer, such that
$$|c| \geq \sqrt{|w|^2 + \chi_k^2(1 + |w|^2)} - \chi_k - \tan^{-1}\left(\frac{
\sqrt{|w|^2 + \chi_k^2(1 + |w|^2)} - \chi_k}{1 + \chi_k \sqrt{|w|^2 + \chi_k^2(1 + |w|^2)}} \right).$$ 
If there exists  $k$, such  that the above inequality is strict, then there are $2k + 1$ geodesics connecting these two points. If $k$ gives the equality, then the number of geodesics is $2k$.
If no such $k$ exists, then there is a unique geodesic.
\end{itemize}
\end{proposition}

\begin{remark}
The number of geodesics in the case (c) above may be difficult to determine, but due to the fact that the value of
$\sqrt{|w|^2 + \chi_k^2 + |w|^2 \chi_k^2} - \chi_k - \tan^{-1}\left(\frac{\sqrt{|w|^2 + \chi_k^2 + |w|^2 \chi_k^2}
- \chi_k}{1 + \chi_k \sqrt{|w|^2 + \chi_k^2 + |w|^2 \chi_k^2}} \right)$, belongs to the interval
 $\Big(\pi k (\sqrt{|w|^2 +1} - 1),\,\, \pi (k + \frac{1}{2}) (\sqrt{|w|^2 +1} - 1)\Big)$,
this number remains between $2k_1 -1$
and $2k_1 +1$, where
$$k_1 = \left\lceil \tfrac{|c|}{\pi(\sqrt{1 + |w|^2} -1)} - \tfrac{1}{2} \right\rceil.$$
\end{remark}

The proof is a long case by case analysis, and therefore, we leave it  to section \ref{proofs}.
Regarding to $\mathrm{(a)}$, we say that two geodesics are {\it geometrically similar},
if one is the image of another under an isometry.
The isometry considered in $\mathrm{(a)}$ is $(c,w) \mapsto (c, e^{i\theta}w)$,
$\theta \in \mathbb{R}$.

Since $\tilde D$ is a strongly bracket generating distribution (i.e. for every vector field $\tilde X_0$, $\tilde D$ and $[\tilde X_0, \tilde D]$ span the entire tangent space at points where $X_0$ is nonvanishing), there
are no abnormal length minimizers.  Every geodesic may be extended indefinitely, therefore the Carnot-Carath\'eodory metric is complete (see \cite{Str1} and \cite{Str2}).
Hence, every point is connected to $\tilde 1$ by a length minimizing geodesic which is normal.
The same holds for $\SU(1,1)$.
The above information along with the proof of Proposition \ref{sRnumbergeo},  leads to the following
result.
\begin{corollary} \label{CorlengthsR}
\begin{itemize}
\item[(a)] If $\tilde g \in e^{\spn \{Z\}}$, then $d(\tilde 1; \tilde g) = |c|$.
\item[(b)] If $\tilde g \in e^{\spn \{X,Y\}}$, then $d(\tilde 1; \tilde g) = \sinh^{-1}|w|$.
\item[(c)] If $\tilde g = (c,w)$ is neither in the vertical, nor in the horizontal space, and $|c| < |w| - \tan^{-1}|w|$,
then
$$d(\tilde{1};\tilde g) = \tfrac{1}{\sqrt{1-\beta_{c,w}^2}} \sinh^{-1}\left(|w| \sqrt{1-\beta_{c,w}^2} \right),$$
where $0<\beta_{c,w}< 1$ is a unique number satisfying
$$|c| = \tfrac{\beta_{c,w}}{\sqrt{1-\beta_{c,w}^2}} \sinh^{-1}\left(|w| \sqrt{1-\beta_{c,w}^2} \right)- \sin^{-1} \left(\tfrac{|w|\beta_{c,w}}{\sqrt{1 + |w|^2}}\right).$$
\item[(d)] If $\tilde g = \left(\pm (|w| - \tan^{-1}|w|) , w \right)$, then $d(\tilde 1;\tilde g) = |w|.$
\item[(e)] If $\tilde g = (c,w)$ satisfies the inequality $ |w| - \tan^{-1} |w|<|c|\leq \frac{\pi}{2}(\sqrt{|w|^2 +1} -1)$, then 
$$d(\tilde{1};\tilde g) = \tfrac{1}{\sqrt{\beta_{c,w}^2-1}} \sin^{-1}\left(|w| \sqrt{\beta_{c,w}^2-1} \right),$$
where $1<\beta_{c,w}\leq \tfrac{\sqrt{1+|w|^2}}{|w|}$ is a unique number satisfying
$$|c| = \tfrac{\beta_{c,w}}{\sqrt{\beta_{c,w}^2-1}} \sin^{-1}\left(|w| \sqrt{\beta_{c,w}^2-1} \right)- \sin^{-1} \left(\tfrac{|w|\beta_{c,w}}{\sqrt{1 + |w|^2}}\right).$$
\item[(f)] If $\tilde g = (c,w)$ satisfies the inequality $\frac{\pi}{2}(\sqrt{|w|^2 +1} -1) < |c|$, then
$$d(\tilde{1};\tilde g) = \tfrac{1}{\sqrt{\beta_{c,w}^2-1}} \left(\pi -\sin^{-1}\left(|w| \sqrt{\beta_{c,w}^2-1} \right)\right),$$
where $1<\beta_{c,w}\leq \tfrac{\sqrt{1+|w|^2}}{|w|}$ is a unique number satisfying
$$|c| = \tfrac{\beta_{c,w}}{\sqrt{\beta_{c,w}^2-1}} \left(\pi- \sin^{-1}\left(|w| \sqrt{\beta_{c,w}^2-1} \right)\right)
- \left(\pi-\sin^{-1} \left(\tfrac{|w|\beta_{c,w}}{\sqrt{1 + |w|^2}}\right)\right).$$
\end{itemize}
\end{corollary}

The details are again left to section \ref{proofs}. 
Notice that in all cases, the distance is independent of the sign of $c$,
and of the argument of $w$.

\begin{remark}
From these results for the universal cover, we easily obtain the following conclusions about the number of geodesics connecting 1 and $g \in \SU(1,1)$. It is
\begin{itemize}
\item[(a)] uncountable (there are countably many geometrically different geodesics), if $g \in e^{\spn \{Z\}}$.
\item[(b)] countable otherwise.
\end{itemize}
\end{remark}

The next result we prove here.
\begin{lemma} \label{lemmashorter}
Let $c_1$ and $c_2$ be two real numbers such that $|c_1| \leq |c_2|$. Then for any $w \in \comp$,
$$d(\tilde 1; (c_1,w) ) \leq d(\tilde 1; (c_2,w)).$$
\end{lemma}

\begin{proof}
If $w = 0$, then this follows from Corollary \ref{CorlengthsR} (a). For $w \neq 0$, observe that the formulas for the distance in Corollary \ref{CorlengthsR} (c) and (e) is increasing with respect to $\beta_{c,w}$, while the distance function in (f) is decreasing with respect to $\beta_{c,w}$. By taking upper and lower limits of the permitted values of $\beta_{c,w}$, we obtain that the cases (b) to (f) are given in the order of increasing distance.

To complete the proof, we need to show that in (c) and (e), $\beta_{c,w}$ is an increasing function with respect to $|c|$, while in (f), $\beta_{c,w}$ is a decreasing relative to $|c|$. We will only show that this holds in (e). The other cases are done similarly.

We need to differentiate
$$f(\beta) = \tfrac{\beta}{\sqrt{1-\beta^2}} \sinh^{-1}\left(|w| \sqrt{1-\beta^2} \right)- \sin^{-1} \left(\tfrac{|w|\beta}{\sqrt{1 + |w|^2}}\right),$$
and see that the derivative is positive. This follows from the computation
$$\tfrac{\partial}{\partial \beta} f(\beta) = \tfrac{1}{1-\beta^2} \left(\tfrac{\sinh^{-1}(|w| \sqrt{1-\beta^2})}{\sqrt{1-\beta^2}} - \tfrac{|w|}{\sqrt{1 + |w|^2(1-\beta^2)}} \right) > 0.$$

\end{proof}

\begin{corollary} \label{ngeoSU}
If $g =(z_1,z_2) \in \SU(1,1)$, then $d(1;(z_1,z_2)) = d(\tilde 1; (\Arg \, z_1, z_2))$. 
\end{corollary}

\begin{proof}
Since being a normal geodesic is a local property, it follows that  any normal geodesic from 1 to $g$ in $\SU(1,1)$, has a unique lift to a curve starting at $\tilde 1$, which is a normal geodesic of equal length to some point in $\pi^{-1}(g)$. This allows us to compute the length by the formula
$$d(1; g) = \min_{\tilde g \in \pi^{-1}(g)} \{d(\tilde 1; \tilde g)\}.$$

To complete the proof of the Corollary, we need to show that for $g = (z_1, z_2)$, the point $(\Arg z_1, z_2) \in \pi^{-1}(g)$ has the minimal distance to $\tilde 1$. Define $c_0 = \Arg z_0$ and let $c_k = c_0 + 2\pi k$, $k \in \mathbb{Z}$. From the formula of the projection \eqref{projwSU11toSU11}, any element in $\pi^{-1}(z_1,z_2)$ is of the form $(c_k, z_2)$. The rest follows from Lemma~\ref{lemmashorter}, since $|c_0| \leq |c_k|$.
\end{proof}

\subsection{The cut and conjugate loci} \label{cutcon}

For  sub-Riemannian Lie groups satisfying Theorem \ref{thmon},
we define a sub-Riemannian analogue of the  exponential map about the identity $1$  by
$$\exp_{sr}:\mathfrak{g} \setminus \mathfrak{k} \rightarrow G,
\qquad X_0 \mapsto e^X_0 e^{- \pr_{\mathfrak{k}} X_0}.$$
For a more general definition of the exponential map in the sub-Riemannian setting,
see \cite{AExp}. We define {\it the conjugate locus} of $G$ from the identity 1, as the set of
critical values of $\exp_{sr}$. We often split the conjugate locus in several sets, defining
the $n$-th conjugate locus by the set $\{ \exp t X_0 \}$, where $X_0 \in \mathfrak{g} \setminus \mathfrak{k}$,
and t is so that there exist exactly $n$ values $0 < t_1 < \cdots < t_n = t$, so that $t_j X_0$ are all critical points.
We define {\it the cut locus} from 1 as the set reachable by
more than one minimizing geodesic.

\begin{corollary}
The cut locus from $\tilde 1$ for $\wSU(1,1)$ is the vertical line. The cut locus from 1 for $\SU(1,1)$ consists of the points where $z_1< 0$, as well as the points apart from the identity satisfying $z_2 = 0$.
\end{corollary}
\begin{proof}
The cut locus for $\wSU(1,1)$ follows from the proof of Proposition \ref{sRnumbergeo} and
Corollary~\ref{CorlengthsR}. For $\SU(1,1)$, if a point $(z_1,z_2)$ is in the cut locus, then either
$c = \Arg \, z_1, w = z_2$ is in the cut locus for $\wSU(1,1)$  (which is the set of points $(c,0)$
for arbitrary $c\neq 0$) or
there exist more than one $\tilde g \in \wSU(1,1)$, such that 
 $\pi(\tilde g) = g$ and $d(\tilde 1 ; \tilde g) = d(1;g)$. From the proof of Corollary \ref{ngeoSU}, this only happens when $\Arg \, z_1 = \pi$, in which there are two points of equal distance.
\end{proof}

The following proposition was proved for $\SU(1,1)$ in \cite{BoRo} (for the isometric case of $\SL(2)$),
but here we generalize it, including the universal cover $\wSU(1,1)$.
\begin{proposition} \label{Propconjugate}
The $n$-th conjugate locus of $\tilde 1 \in \wSU(1,1)$ consists of the vertical line, if $n$ is odd.
If $n = 2j$, then it consists of the points given by the equation
$$|c| = \sqrt{|w|^2 + \chi_j^2 + |w|^2 \chi_j^2} - \chi_j - \tan^{-1}\left(\tfrac{\sqrt{|w|^2 + \chi_j^2 + |w|^2 \chi_j^2}
- \chi_j}{1 + \chi_j \sqrt{|w|^2 + \chi_j^2 + |w|^2 \chi_j^2}} \right).$$
\end{proposition}
\begin{proof}
First, observe that from the definition of the exponential map \eqref{wSU11exp} and Remark \ref{chi} in Section \ref{notationsec}, $d\exp_{sr}$ exists only for  $\check a \neq 0$. 
Put $w = u + iv$. Then we have
$$(c,u+iv) = \exp_{sr}\left(r\cos \theta X + r \sin \theta Y + a_3 Z\right)$$
$$= \left(-\scrT_a + a_3, - r \scrS_{r^2 - a_3^2}\cos(\theta + a_3) +i r \scrS_{r^2 - a_3^2} \sin(\theta + a_3)\right).$$
The values of the elements of the matrix
$$\left(\begin{array}{ccccc}
\frac{\partial}{\partial \theta} c &\,& \frac{\partial}{\partial \theta} u &\,& \frac{\partial}{\partial \theta} v \\
\frac{\partial}{\partial r} c &\,& \frac{\partial}{\partial r} u &\,& \frac{\partial}{\partial r} v \\
\frac{\partial}{\partial a_3} c &\,& \frac{\partial}{\partial a_3} u &\,& \frac{\partial}{\partial a_3} v \\
\end{array} \right),$$
become
$$\left(\begin{array}{ccccc}
0 &\, &r S \sin(\theta + a_3) &\,&r S \cos(\theta + a_3) \\
- \frac{a_3 r}{r^2 - a_3^2} \frac{1-CS}{C^2 + a_3 S^2} &\,& -\frac{(r^2 C- a_3^2S)\cos(\theta + a_3)}{r^2 - a_3^2} &\,&
\frac{(r^2 C- a_3^2S)\sin(\theta + a_3)}{r^2 - a_3^2} \\
1 + \frac{a_3^2 - r^2 CS}{(r^2 - a_3^2)(C^2 + a_3^2 S^2)} &\,& 
r S \sin(\theta) + \frac{ra_3(C-S)\cos(\theta + a_3)}{r^2 - a_3^2} &\,&
r S \cos(\theta) - \frac{ra_3(C-S)\sin(\theta + a_3)}{r^2 - a_3^2} \\
\end{array}\right).$$
Here, we have simplified by writing just $S= \scrS_{r^2 -a_3^2}$ and $C = \scrC_{r^2 -a_3^2}$.
The determinant of the above matrix is $\frac{r^3 S (C-S)}{r^2 - a_3^2}$.
The value of $S$ vanishes at the points, for which  $a_3^2 - r^2 = \pi^2 k^2, k \in \mathbb{N}$,
and the image of such points is the vertical line
(see the proof for Proposition~\ref{sRnumbergeo} for more details).
Moreover, $C-S$ vanishes only for $- \check a = \chi_j^2$.
Hence, for a generic $X \in \mathfrak{g} \setminus \mathfrak{k}$, the point $tX$ can 
be  singular for $\exp_{sr}$, only if $\check a < 0$. Let us use the normalization 
$\check a = -\pi^2$.
Let  the $n$-th value $t_n$ be such that $t_nX$, is a singular point.
Then it is clear that
$$t_{2j-1} = j, \qquad t_{2j} = \frac{\chi_j}{\pi}, \qquad j \in \mathbb{N}. $$
Taking the image of these values under the map $\exp_{sr}$, we have the result.
\end{proof}

\begin{remark}
Let $(c,w) \in \wSU(1,1)$ be a point that does not belong to the vertical line (i.e. $w \neq 0$).
Notice that if  the value of $|c|$ from Proposition~\ref{sRnumbergeo} increases (or the value of $|w|$ decreases), then (and only then) the number of geodesics increases when we pass through the even-indexed conjugate loci
(see Figure 2).
\end{remark}
\begin{figure}[ht!]
\centering
\scalebox{0.60}{\includegraphics{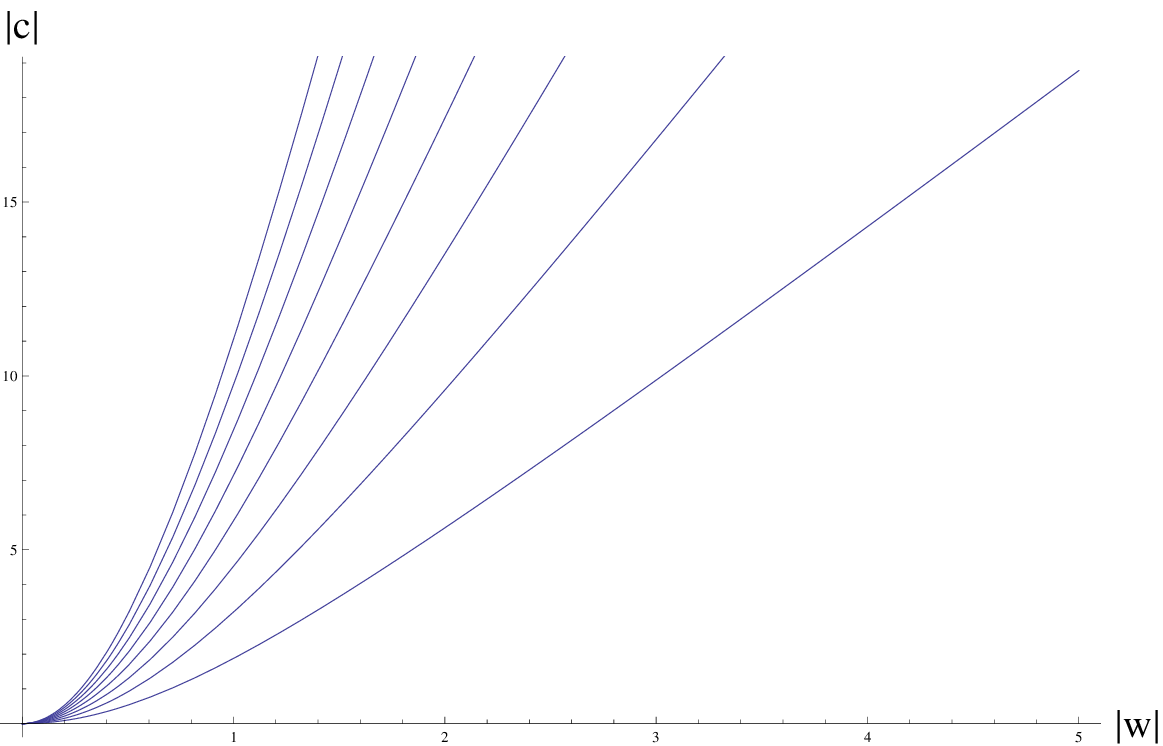}}
\caption[]{The above figure shows the eight first even conjugate loci. The horizontal axis corresponds
to the values of $|w|$, and the values of $|c|$ occupy the vertical axis.
Observe that a point under all the graphs is reachable by a unique geodesic. Then there are two geodesics
in the 2-nd conjugate locus (that is the first graph).
Between the 2-nd and the 4-th conjugate locus there are three, and so on.}
\end{figure}

Projecting the result in Corollary \ref{Propconjugate} by \eqref{projwSU11toSU11}.
\begin{corollary}
The $n$-th conjugate locus of $1 \in \SU(1,1)$ for odd $n$  is the vertical line, and
if $n = 2j$, then it consists of the points
$$z_1 = \frac{(-1)^j \left(1 \mp i \sqrt{|z_2|^2 (1 + \chi_j^2) + \chi_j^2}\right)
e^{\pm i\sqrt{|z_2|^2 (1 + \chi_j^2) + \chi_j^2}}}{\sqrt{1 + \chi_j^2}}.$$
\end{corollary}

\section{Sub-Lorentzian geometry on $\wSU(1,1)$}
\subsection{Sub-Lorentzian maximizers and geodesics on $\wSU(1,1)$}
In contrast with  the sub-Riemannian case, we only know that the  relative maximizers exist  locally.
We have no guarantee of global existence of maximizers.
Let us consider the distribution $\widetilde E = \mathrm{span}\{\widetilde Y, \widetilde Z \}$,
with the metric $\tilde \rho$
restricted to $\widetilde E$.
We formulate an optimal control problem of maximizing
$$J(u,\tilde g_u) = \int_0^{\tau} \sqrt{u_1^2 - u_2^2} \, dt,
\qquad u \in U = \left\{u \in \mathbb{R}^2 | \, u_1 = \sqrt{1 + u_2^2} \right\},$$
where $\dot{\tilde g}_u = f(u,\tilde g) = u_1 \widetilde Z(\tilde g) + u_2 \widetilde Y(\tilde g)$.

We only do the case when $\varkappa =1$. If $\varkappa = 0$ one gets the same results,
so it follows that there are no strictly abnormal geodesics.
In order to use PMP, let us define the pseudo-Hamiltonian
$$\mathscr{H}(u, p_1X^* + p_2 Y^* + p_3 Z^*) = \sqrt{u_1^2 - u_2^2} + u_1 p_3 + u_2 p_2$$
$$= 1 + p_3 \sqrt{1 + u_2^2} + p_2 u_2.$$
where $u \in U = \{(u_1,u_2) \in \mathbb{R}^2 | \, u_1 = \sqrt{1 + u_2^2}\}$.
For the existence of $H(p) = \max_{u \in U} \mathscr{H}(p,u)$,
we need that $|p_3| > |p_2|$ and that $p_3< 0$. In this case
$$H(p_1X^* + p_2 Y^* + p_3 Z^*) = 1 - \sqrt{p_3^2 -p_2^2}.$$
This is from the fact that the optimal control is
$$u^*  = \left( \tfrac{-p_3}{\sqrt{p_3^2 - p_2^2}} , \tfrac{p_2}{\sqrt{p_3^2 - p_2^2}} \right).$$
Since $H$ cannot be extended to $T^*M \setminus s_0(M)$ (the tangent bundle with the zero-section removed), we solve the equation using the pseudo-Hamiltonian function $ \mathscr{H}$ instead,
$$\dot{p}_1 = \{p_1, \mathscr{H}(u,\lambda) \} |_{u = u^*} =
-2 (u_1 p_2 + u_2 p_3) |_{u = u^*} = 0,$$
$$\dot{p}_2 = \{p_2, \mathscr{H}(u,\lambda) \} |_{u = u^*} =
2 u_1 p_1 |_{u = u^*} = - \tfrac{2 p_1 p_3}{\sqrt{p_3^2 - p_2^2}},$$
$$\dot{p}_3 = \{p_3, \mathscr{H}(u,\lambda) \} |_{u = u^*} =
 - 2 u_2 p_1 |_{u = u^*} = - \tfrac{2 p_1 p_2}{\sqrt{p_3^2 - p_2^2}},$$
$p_3^2 - p_2^2$ is a first integral in this case, and from the condition
 $\mathscr{H}(u, p_1 X^* + p_2 Y^* + p_3 Z^*) = 0$, $p_3^2 - p_2^2$ is equal to 1.
 We find the solutions (denote $p_{jo} := p_j(0)$)
 $$p_1(t) \equiv p_{1o},$$
 $$p_2(t) =  p_{2o} \cosh\left( 2p_{1o} t \right) 
-  p_{3o} \sinh\left( 2p_{1o} t \right),$$
$$p_3(t) =  p_{3o} \cosh\left( 2p_{1o}t \right) 
-  p_{2o} \sinh\left( 2p_{1o} t \right).$$
Observe that $|p_{3o}| > |p_{2o}|$, and  $p_{3o} < 0$. Note also that $p_{3o}$ has to be equal to
$- \sqrt{1 + p_{2o}^2}.$
In order to simplify our calculation, we first solve it for $\SU(1,1)$, and then lift it to $\wSU(1,1)$.
We have to solve the following differential equation.
\begin{equation} \label{hameq}
\dot{g} = \left. \tfrac{\partial \mathscr{H}}{\partial p} \right|_{u = u^*}
= - p_{3}Z(g) + p_{2}Y(g),
\end{equation}
Make the following observation that
$$e^{a_1 X} = \left( \begin{array}{cc} \cosh a_1 & - \sinh a_1 \\ - \sinh a_1 & \cosh a_1
\end{array} \right).$$
From this we notice that
$$e^{- a_1 X} Z = Z e^{a_1 X} = \cosh a_1 Z + \sinh a_1 Y,$$
$$e^{- a_1 X} Y = Y e^{a_1 X} = \cosh a_1 Y + \sinh a_1 Z,$$
$$e^{a_1 X} X = X e^{a_1 X} = \cosh a_1 X + \sinh a_1 \cdot 1.$$
If we expand equation \eqref{hameq}, we can write
$$\dot{g} = -\left(p_{3o} \cosh\left( 2p_{1o}t \right)  -  p_{2o} \sinh\left( 2p_{1o}t \right)\right) Z(g)
+ \left(p_{2o} \cosh\left( 2p_{1o} t \right) -  p_{3o} \sinh\left( 2p_{1o} t \right)\right) Y(g)$$
$$= g( - p_{3o} Z e^{2 t p_{1o} X} + p_{2o} Y e^{2 t p_{1o} X} )
= g e^{- t p_1 X}( - p_{3o} Z + p_{2o} Y)e^{t p_{1o}X},$$
and multiplying from the right by $e^{- t a_1 X}$ and adding
$ - p_{1o} g X e^{- tp_{1o} X} =  - p_{1o} g e^{- tp_{1o} X}X$ on both sides,
we obtain
$$\dot{g} e^{- t p_{1o} X} - p_{1o} g X e^{- t p_{1o} X}
= \tfrac{\partial}{\partial t} \left({g} e^{- t p_{1o} X} \right)
= g e^{- t p_{1o} X} (- p_{3o} Z + p_{2o} Y - p_{1o} X).$$
It follows that
$$g e^{- t p_{1o} X} = e^{t(-p_{1o}X + p_{2o} Y - p_{3o} Z)},$$
i.e.,
$$g  = e^{t(-p_{1o}X + p_{2o} Y - p_{3o} Z)} e^{t p_{1o} X}.$$

If we lift this curve and consider Theorem \ref{thmon}, then we easily get to the following proposition.
\begin{proposition}
Assume that there exists a global length maximizing curve between $\tilde 1$ and $\tilde g$. Then,
this curve is a timelike, future-directed normal geodesic on the form
$$e^{t(a_1 X + a_2 Y + a_3 Z)} e^{- t a_1 X},$$
where $a_3 > 0$ and $|a_2| < a_3$.
\end{proposition}
Explicitly, by \eqref{wSU11exp}, these geodesics have the form
\begin{eqnarray*}c(t) = &-& \scrT_{t a} + \tan^{-1} \left( \tfrac{a_2 \scrC_{\check a}(t) \scrS_{\check a}(t) - a_1 a_3 \scrS_{\check a}(t)^2}
{(1 + |a_1 + i a_2|^2 |\scrS_{\check a}(t)|^2) \coth(a_1 t) - a_1 \scrC_{\check a}(t) \scrS_{\check a}(t) - a_2 a_3 \scrS_{\check a}(t)^2} \right),
\end{eqnarray*}
$$w(t) = (\scrC_{\check a}(t) - i a_3 \scrS_{\check a}(t)) \sinh(a_1 t) - (a_1 - ia_2) \scrS_{\check a}(t)\cosh(a_1t).$$
The projection of them to $\SU(1,1)$ is given by
$$z_1 = (\scrC_{\check a}(t) - i a_3 \scrS_{\check a}(t)) \cosh(a_1 t) - (a_1 - ia_2) \scrS_{\check a}(t)\sinh(a_1t),$$
$$z_2 = (\scrC_{\check a}(t) - i a_3 \scrS_{\check a}(t)) \sinh(a_1 t) - (a_1 - ia_2) \scrS_{\check a}(t)\cosh(a_1t).$$

We will discuss the situation when these curves correspond to length maximizers in Section~\ref{future}.

\subsection{Number of geodesics}
The result in the sub-Lorentzian case is  more complicated, than in the sub-Riemannian case.
So we have to give some definitions in order to describe geodesics in a reasonable way. Write $\mathbb N_0$ for the collection of nonnegative integers.
First, let us define a function $f_\pm: \mathbb{R} \times \mathbb{N}_0 \rightarrow \mathbb{R}$, by
$$f_\pm(s; k) = \left\{ \begin{array}{ll}
\sqrt{1 - s^2} & \text{if } |s| < 1, \\
\mp \sinh\left(\frac{(2k \pm 1) \pi}{2} \sqrt{s^2 -1} \right) &
\text{if } |s| \geq 1. \end{array} \right.$$
Also, we define the numbers $\omega_k, k \in \mathbb{N}$ as the numbers satisfying the equations
$$\frac{\omega_k}{\sqrt{1 -\omega_k^2}} - \sin^{-1} \omega_k = \pi k, \qquad 0 < \omega_k < 1.$$
Finally, we define the function $$F(s, \omega) = \sqrt{1 -\omega^2} \cosh
\sqrt{\frac{s^2 - \omega^2}{1 - \omega^2}} - \sqrt{s^2 - \omega^2}
\sinh \sqrt{\frac{s^2 - \omega^2}{1 - \omega^2}}.$$

Let us construct the following subsets of $\wSU(1,1)$:
$\Omega_0$ consists of all points
$$c = \Arg\left(x_1 - i \sqrt{s^2 + y_2^2}\right), \qquad w = \pm \sqrt{x_1^2 + s^2-1} + i y_2,$$
where $y_2 \in \mathbb{R}$, $s >0$, and $x_1 \in (-\infty,\, -\sqrt{1-s^2}]\bigcup [\sqrt{1-s^2}, \,
1 - \frac{s^2}{2})$ when $s < 1$, and $x_1 \in  (-\infty,\, 1 - \frac{s^2}{2})$ otherwise.
Further, for $k \in \mathbb{N}$, define $\Omega_k$ as a set of  all points
$$c = - \pi k + \Arg\left(x_1 - i \sqrt{s^2 + y_2^2}\right), \qquad w = \pm \sqrt{x_1^2 + s^2-1} + i y_2,$$
where $y_2 \in \mathbb{R}$, $s >0$, and $x_1 \in \mathbb{R} \setminus \left(-\sqrt{1-s^2},
\sqrt{1 - s^2}\right)$ for $s < 1$, and $x_1 \in (-\infty, F(s,\omega_k)]$ otherwise.
We split them into the following sets
\begin{itemize}
\item $A_k$ consist of all points
$$c = - \pi k + \Arg\left(x_1 - i \sqrt{s^2 + y_2^2}\right), \qquad w = \pm \sqrt{x_1^2 + s^2-1} + i y_2,$$
where $y_2 \in \mathbb{R}$, $s >0$, and 
\begin{itemize}
\item $x_1 \in \mathbb{R} \setminus \left(-\sqrt{1-s^2},\sqrt{1 - s^2}\right)$, when $s \leq \omega_k$,
\item $x_1 \in (-\infty, -\sqrt{1 - s^2}] \cup (F(s,\omega_k), \infty)$, when $\omega_k < s < 1$,
\item $x_1 \in (- \infty, 0] \cup \{F(1,\omega_k)\}$, when $s =1$,
\item $x_1 \in \{F(s,\omega_k)\}$, when $s> 1$.
\end{itemize}

\item $B_k$ consist of all points
$$c = - \pi k + \Arg\left(x_1 - i \sqrt{s^2 + y_2^2}\right), \qquad w = \pm \sqrt{x_1^2 + s^2-1} + i y_2,$$
where $y_2 \in \mathbb{R}$, $s >0$, and 
\begin{itemize}
\item $x_1 \in \{\sqrt{1 -s^2}, F(s,\omega_k)\}$, when $s < \omega_k $,
\item $x_1 \in (0, F(s,\omega_k))$, when $s =1$,
\item $x_1 \in (- \infty, F(s;k))$, when $s > 1$.
\end{itemize}

\item $C_k$ consist of all points
$$c = - \pi k + \Arg\left(x_1 - i \sqrt{s^2 + y_2^2}\right), \qquad w = \pm \sqrt{x_1^2 + s^2-1} + i y_2,$$
where $y_2 \in \mathbb{R}$, $\omega_k < s < 1$, and $x_2 \in (\sqrt{1 - s^2}, F(s,\omega_k))$.
\end{itemize}
Let $\Xi_k = \{ (-\pi k, x_2) \in \wSU(1,1) | \, x_2 \in \mathbb{R} \}$.
Finally, define
$$\Omega := \Omega_0 \, \cup \, \bigcup_{k= 1}^\infty (\Omega_k \, \cup \, \Xi_k),
\qquad \Xi := \bigcup_{k = 1}^\infty \Xi_k,
$$
$$
A := \Omega_0 \, \cup \, \bigcup_{k = 1}^\infty A_k, \qquad B := \bigcup_{k =1}^\infty B_k, 
\qquad C := \bigcup_{k = 1}^\infty C_k.$$

\begin{proposition} \label{maxgeo}
There are timelike, future directed geodesics connecting  $\tilde 1$ and $\tilde g$,
only if  $\tilde g \in \Omega$. More precisely,
\begin{itemize}
\item the geodesic is unique if $\tilde g \in A$;
\item there are two geodesics if $\tilde g \in B$;
\item there are three geodesic if $\tilde g \in C$;
\item there are countable many geodesics if $\tilde g \in \Xi$.
\end{itemize}
\end{proposition}
The proof for this is a case by case analysis, which we leave for section \ref{proofmaxgeo} .
\begin{remark}
Let $\chi_k$ be the constants defined in \eqref{defchik} in Section \ref{notationsec}.
Then $\omega_k = |\sin \chi_k|$, so these numbers play a role here as well.
However, unlike what was the case when we considered the sub-Riemannian geodesics, the mapping
$$a_1 X + a_2 Y + a_3 Z \mapsto e^{a_1 X + a_2 Y + a_3 Z} e^{-a_1X},$$
does not have any critical points satisfying $\check a = - \chi_k^2$.
\end{remark}

We also remark the following interesting comparison with Lorentzian geometry.
\begin{proposition}\label{Prop6}
The set reachable by Lorentzian geodesics starting at $\tilde 1$, is the planes
$\{(c,w) \in \wSU(1,1) | \, c = - \pi k -\frac{\pi}{2}, k \in \mathbb{N} \}$, along with all
points satisfying
$c< 0$, $|w|<\tan |c|$, and $(c,w) = (- \pi k, 0), k \in \mathbb{N}$.
More precisely,
\begin{itemize}
\item[(a)] there are uncountably many geodesics if $w = 0$
(all geometrically similar, except for the case $c = 0 \text{ mod } \pi$, when they are 
uncountably many geometrically
different).
\item[(b)] there is a unique geodesic if $w \neq 0$.
\end{itemize}
\end{proposition}

\begin{proof}
A time-like geodesic $\tilde g(t)$, starting at $\tilde 1$ satisfies $\tilde{g}(1) = e^{a_1X + 
a_2 Y + a_3Z}$,
and we know that $a_3 > 0$ and $a_3^2 > a_1^2 + a_2^2$ (and hence $\check a < 0$).
We rewrite the geodesics in terms of the parameters $k,\alpha' , r$ and $\theta$, where
$\sqrt{-\check a} = \pi k + \alpha', k \in \mathbb{N}, \alpha' \in (- \frac{\pi}{2}, 
\frac{\pi}{2}]$ and
$a_1 + i a_2 = r e^{i\theta}$. The geodesics in these coordinates look like
$$c = - \tan^{-1} \left(\frac{\sqrt{r^2 + (\pi k + \alpha')^2}}{\pi k + \alpha'} \tan 
\alpha' \right)
- k \pi, \qquad w = - \frac{r}{\pi k + \alpha'} e^{i\theta} (-1)^k \sin \alpha' .$$
Now $\theta$ can be determined by $\arg w$, and $k = \lceil - \frac{c}{\pi} -\frac{1}{2} 
\rceil$.
Rewriting the equation for $|w|$ as
$$|w|^2 + \sin^2 \alpha' = \frac{r^2 + (\pi k + \alpha')^2}{(\pi k + \alpha')^2} \sin^2 
\alpha',$$
and inserting it in the equation for $c$, we get that
$$\alpha' = - \sin^{-1}\left( \tan(c) \sqrt{1 - \frac{|w|^2}{\tan^2 c}} \right),$$
for $c\neq  \frac{\pi n}{2}$.
For the remaining cases, $\alpha' = 0$ for $c = \pi n$, and $\alpha ' = 
\frac{\pi}{2}$
for $c = \frac{\pi}{2} +\pi n$. It is obvious that there are  solutions only if 
$|\tan c| \geq |w|$.
The equality can never be attained, because $\sgn \tan c = - \sgn \alpha'$, which follows from
$$\tan c = - \frac{\sqrt{r^2 + (\pi k + \alpha')^2}}{\pi k + \alpha'} \tan \alpha' .$$
The equation
$$r^2 = \frac{(\pi k + \alpha')^2 |w|}{\sin^2 \alpha'},$$
determines $r$ for $\alpha' \neq 0$. When $\sin \alpha' = 0$, there are no restrictions 
on $r$.
Also, if $|w| = 0$, there are no restrictions on $\theta$.
\end{proof}

Proposition~\ref{Prop6} yields that the set of points reachable by the Lorentzian geodesics starting at $\tilde 1$ 
neither is contained,
nor contains the set reachable from $\tilde 1$ by the sub-Lorentzian geodesics.
This contrasts  the fact that the Lorentzian timelike future always contains the 
sub-Lorentzian one.

\subsection{Lorentzian and sub-Lorentzian timelike future}\label{future}
Using the information we have collected, concerning sets reachable by geodesics, we obtain some results for the timelike future in general.

\vspace{0pt}
\begin{figure}
\begin{center}
\begin{pspicture}(0,0)(7,7)
\put(-3,-2){\epsfig{file=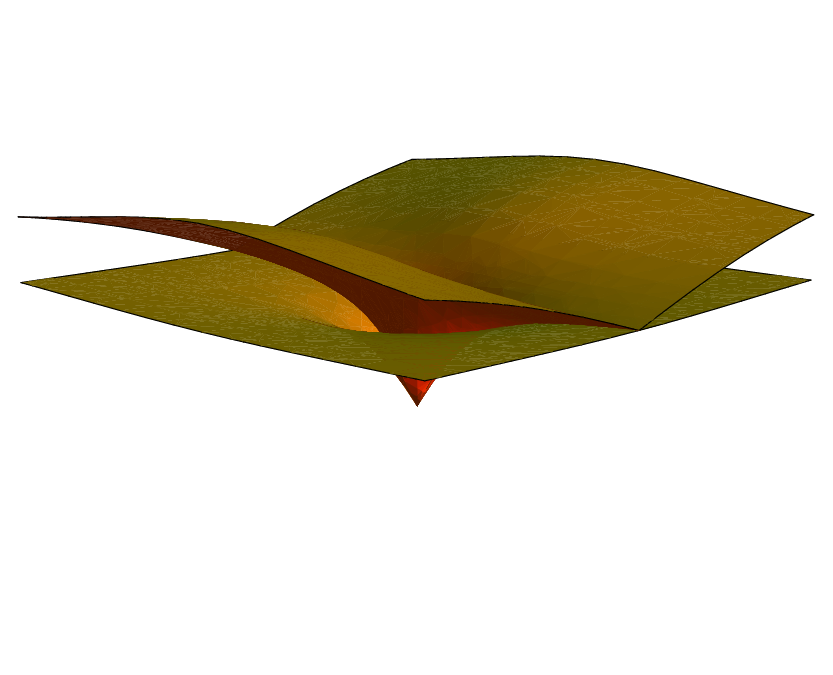,height=9cm}}
\rput(2.7,6){$c$}
\rput(8,2.5){$\re w$}
\rput(5,0.6){$\im w$}
\rput(2.5,1.3){$0$}
\psline[linewidth=0.3mm]{->}(2.5,3.4)(2.5,6)
\psline[linewidth=0.3mm]{->}(-1,0.63)(8,3.03)
\psline[linewidth=0.3mm]{->}(0,2.2)(5,1)
\end{pspicture}
\end{center}
\caption[]{The above figure shows the lower bounds of Lorentzian and sub-Lorentzian timelike future for $\wSU(1,1)$,
showing the sub-Lorentzian future contained in the Lorentzian one.
The negative direction of $c$ in the picture is upward.
Remark that along the
imaginary $w$-axis, $\mathcal{I}^+(\tilde 1, \tilde \rho|_{\tilde E})$ is tangent to $\mathcal{I}^+(\tilde 1, \rho)$.}
\end{figure}
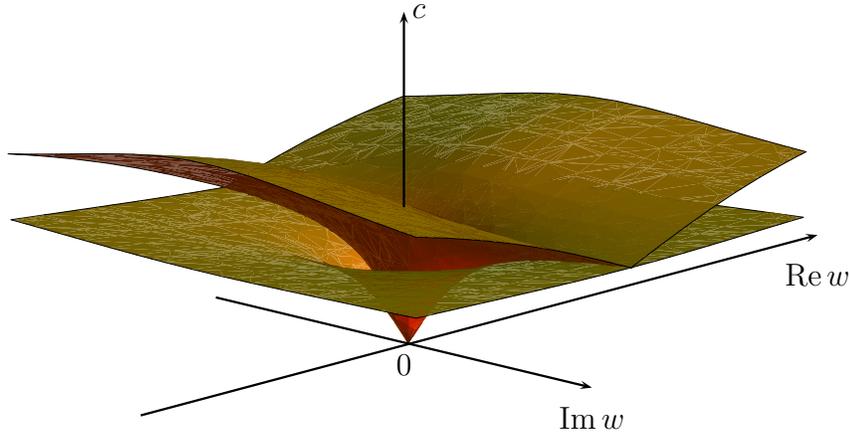

\begin{proposition}
\begin{itemize}
\item[(a)] The timelike future of $\wSU(1,1)$ with respect to the Lorentzian metric $\tilde \rho$,
is given by
$$\mathcal{I}^+(\tilde 1, \tilde \rho) = \left\{(c,w) \in \wSU(1,1)\, | \, c < - \tan^{-1} 
|w| \right\}.$$
\item[(b)] The timelike future of $\wSU(1,1)$ with respect to the sub-Lorentzian metric $\tilde 
\rho|_{\tilde{E}}$,
is given by
$$\mathcal{I}^+(\tilde 1, \tilde \rho|_{\tilde E}) =
\{(c,w) \in \wSU(1,1)\, | \, c < \Arg(1- |\re w| - i \sqrt{(\im w)^2 + 2 |\re w|}) \}.$$
\end{itemize}
\end{proposition}
A picture of this is shown in Figure 3.
\begin{proof}
Let us denote  by $Q = \left\{(c,w) \in \wSU(1,1)| \, c < - \tan^{-1} |w| \right\}$, and 
$$\widehat{Q} =
\{(c,w) \in \wSU(1,1)| \, c < \Arg(1- |\re w| - i \sqrt{(\im w)^2 + 2 |\re w|}) \}.$$
In both  sub-Lorentzian and Lorentzian settings, information about the geodesics, imply that $Q$ and $\widehat{Q} $ are included in their respective timelike futures.
To show the opposite inclusion, remark that in both  Lorentzian and sub-Lorentzian 
settings,
we use the fact that there exist relative length maximizers locally.
So  there exists a neighborhood $V$ of $\tilde 1$, such that
$$\mathcal{I}^+(\tilde 1, \tilde \rho|_V) = V \cap Q.$$
Define $U = V \cap Q$. Now
from left invariance, we know that $\mathcal{I}^{+}(\tilde g, \tilde \rho) = L_{\tilde g}
\mathcal{I}^{+}(\tilde 1, \tilde \rho)$, and it follows that
$$\mathcal{I}^{+}(\tilde 1, \tilde \rho) = \bigcup_{\tilde g \in \mathcal{I}^+(\tilde 1, 
\tilde \rho)} L_{\tilde g} U.$$
Hence it is sufficient to show that $L_{\tilde g} U \subseteq Q$, for every
$\tilde g \in Q$. Since every timelike curve from $\tilde 1$ to a point outside $Q$ has 
to pass through the boundary
$\partial Q$, it follows by continuity that it is sufficient to show that $L_{\tilde g_0} 
\overline{U} 
\subseteq L_{g_0} \overline{Q} \subseteq \overline{Q}$, for every $\tilde g_0 \in 
\partial Q$. Here, $\overline{U}$ means the closure of $U$.
Finally, since $L_{\tilde g_0}$ is an isomorphism, we have that $\partial L_{\tilde g_0} 
\overline{Q} = L_{\tilde g_0}
\partial Q$, so all  remaining arguments turn down to show  that $\tilde g_0 \tilde g \in 
\overline{Q}$ for
every $\tilde g_0, \tilde g \in \partial Q$. The same holds for $\widehat{Q}$.

In order to prove (a) we show that there are  timelike geodesics connecting every point satisfying $- 
\frac{\pi}{2} \leq c < - \tan |w|$ with the origin.
Let $(c,w)$ be an arbitrary point satisfying $c < - \frac{\pi}{2}$. Then we can construct a 
timelike curve from $\tilde 1$,
by taking a geodesic from $\tilde 1$ to $(\frac{\pi}{2}, i w e^{-ic})$, and then continue 
by a left translation
of a geodesic from $\tilde 1$ to $(c + \frac{\pi}{2},0)$.

Let $\tilde g_0, \tilde g \in \partial Q$, $\tilde g_0 = (- \tan^{-1} r_0 , r_0 
e^{i\theta_0}),
\tilde g = (- \tan^{-1} r , r e^{i\theta})$.
First, if $r \geq 1$ and $r_0 \geq 1$, then
$$\left(- \tan^{-1} r_0, r_0 e^{i \theta_0}, - \tan^{-1} r , r e^{i \theta} \right) = (C, 
W),$$
$$C = - \tan^{-1} r_0 - \tan^{-1} r  - \tan^{-1} \left(\frac{(r r_0 -1) \sin\vartheta
+ (r + r_0) \cos \theta}{\frac{1}{r r_0}+ \frac{r}{r_0} + \frac{r_0}{r} + r r_0 +
(r + r_0) \sin\vartheta - (r r_0 -1) \cos\vartheta } \right), $$
$C\leq - \frac{\pi}{2},$
so in this case $(C,W) \in Q$. Here $\vartheta = \theta_0 - \theta$.

We now turn to the case when either $r_0$ or $r$ is less than 1. If we denote
$$(x_1 + i y_1, z_2) := \pi(\tilde g_0) \pi(\tilde g),$$
then we have
$$x_1 = 1 - r r_0 + r r_0 \cos \vartheta,$$
$$y_1 = - r_0 - r + r r_0 \sin \vartheta,$$
$$|z_2|^2 = r^2 + r_0^2 + 2 r^2 r_0^2 + 2r r_0 \left( (1 - r r_0) \cos \vartheta - (r + 
r_0) \sin \vartheta \right)$$
$$= 2 ( 1- r r_0) x_1 + 2 (r + r_0) |y_1| - (r + r_0)^2 - 2 (1 - r r_0).$$
Since  $r < 1$ or $r_0 < 1$, we conclude that $y_1 < 0$.
It follows that $\Arg( x_1 + i y_1) \leq - \tan^{-1} |z_2| = \Arg(1 -i |z_2|)$, because
$x_1 \leq 1$.
As a consequence, $g_0 g \in \overline{Q}$.

Turning to (b), let us observe that $\Omega_0$ consist of all points, satisfying
$$-\pi +\tan^{-1}\left(\sqrt{\frac{(\im w)^2}{1 + (\re w)^2}}\right) <c < \Arg(1- |\re w| - i 
\sqrt{(\im w)^2 + 2 |\re w|}).$$
Clearly all the points in $\Omega_0$ are in the timelike future.
Pick up any point $(c,w)$ satisfying $c < \Arg(1- |\re w| - i \sqrt{(\im w)^2 + 2 |\re w|})$.
Let $\gamma$ be a left translation by $(0, w e^{-ic})$ of any geodesic connecting $\tilde 1$ and 
$(c, 0)$.
The endpoint of $\gamma$ is  $(c,w)$. Since $c < \Arg(1- |\re w| - i \sqrt{(\im w)^2 
+ 2 |\re w|})$, the curve
$\gamma$ must at some point intersect $\Omega_0$. Pick up any point $\tilde g_1 \in \gamma 
\cap \Omega_0$,
travel from $\tilde 1$ to $\tilde g_1$  along a geodesic, and continue along $\gamma$ to 
$(c,w)$.

Let $\tilde g_0, \tilde g \in \widehat{Q}$ be of the form
$$g_0 = \left(- \Arg(1- |x| - i \sqrt{y^2 + 2|x|}), x + i y \right), \,\,\, \tilde
g = \left(- \Arg(1- |u| - i \sqrt{v^2 + 2|u|}), u + i v \right).$$
Assume first that $|u| \leq \frac{1}{2}$.
If we again denote
$\pi(\tilde g_0) \pi(\tilde g) = (x_1 + i y_1, x_2 + i y_2)$, we then obtain 
$$x_1 = (1 - |x|)(1 - |u|) + xu + vy - \sqrt{y^2 + 2 |x|} \sqrt{v^2 + 2|u|},$$
$$x_2 = (1 - |x|) u + v \sqrt{y^2 + 2 |x|} + x (1 - |u|) - y \sqrt{v^2 + 2|u|}.$$
We only need to show that $1 - x_1 - |x_2| \geq 0$. Denote $\sigma = \sgn x_2$.
Then
$$1 - x_1 - |x_2| = (|x| - \sigma x) + (|u| - \sigma u) - (|u| - \sigma u) (|x| - \sigma 
x)
+ (\sqrt{v^2 + 2 |u|} - \sigma v) (\sqrt{y^2 + 2 |x|} + \sigma y) $$
$$\geq (|x| - \sigma x) + (|u| - \sigma u) - (|u| - \sigma u) (|x| - \sigma x) \geq 0,$$
since $0 \leq (|u| - \sigma u) \leq 1$. If $|u| > \frac{1}{2}$, the equality
$$\left(- 1- |u| - i \sqrt{v^2 + 2|u|}, u + i v \right)$$
$$= \left(1 - \frac{i}{2} \frac{1}{\sqrt{v^2 + 2|u|} + \sigma_1v}, - \frac{i}{2} 
\frac{\sigma_1}{\sqrt{v^2 + 2|u|} + \sigma_1v} \right)(1 - |u_1| - i \sqrt{v_1^2 + 
2 |u_1|}, u_1 + i v_1),$$
where $\sigma_1 = \sgn u$, and
$$u_1 = u - \frac{\sigma_1}{2}, \qquad v_1 = \frac{\sigma_1}{2} \left(
\sqrt{v^2 + 2|u|} + \sigma_1 - \frac{2|u|-1}{\sqrt{v^2 + 2|u|} - \sigma_1 v}\right),$$
tells us that $\tilde g$ may be written as a finite product of boundary elements, all 
with real part of
the second coordinate of absolute value less than or equal to $\frac{1}{2}$.
\end{proof}

Since we know that locally, the Lorentzian or sub-Lorentzian causal future is the closure of the timelike 
future (see \cite{Groch3}), it follows from
left invariance that $\mathcal{J}^+(\tilde 1, \tilde\rho)$ and $\mathcal{J}^+(\tilde 1, \tilde\rho|_{\tilde E})$ are the closure of 
$\mathcal{I}^+(\tilde 1, \tilde \rho)$ and $\mathcal{I}^+(\tilde 1, \tilde \rho|_{\tilde E})$, respectively.

\begin{lemma}
Both with respect to the Lorentzian $\tilde \rho$ and the sub-Lorentzian metric $\tilde 
\rho|_{\tilde E}$,
$\wSU(1,1)$ is strongly causal.
\end{lemma}

\begin{proof}
The proof is the same for both the Lorentzian  and sub-Lorentzian cases. By left 
invariance, it is sufficient to find a strongly causal neighborhood.
Let $U_{c_0}, 0 < c_0 < \frac{\pi}{2}$ be the set of all points $(c,w)$ in the timelike 
future of $\tilde 1$ satisfying
$c > - c_0$. Observe that inequality $|w| < \tan |c|$ must hold at the same time.
Since, the causal future of all the elements in $U_{c_0}$ is contained in the causal future of $\tilde 1$,
any causal curve must exit through the surface $\{(- c_0, w_0)| \, |w_0| \leq \tan^{-1} c_0\}$.
Then we need to show that for any such element 
$L_{(c_0, w_0)} U_{c_0} \cap U_{c_0} = \emptyset$. Observe that if $(c,w) \in U_{c_0}$ and
$(-c_0, w_0) (c,w) = (C,W)$, then
$$C = -c_0 + c + \Arg(\sqrt{(1 + |w_0|^2)(1+|w|^2)} + w_0 w e^{-i(c + c_0)})$$
$$\leq - c_0 + c  + \tan^{-1} \left(\frac{|w_0 w|}{\sqrt{(1+ |w_0|^2)(1+|w|^2)}} \right) <
- c_0 + c + \tan^{-1} (\sin c_0 \sin |c|)$$
$$ < - c_0 + c + \tan^{-1} (\sin |c|) < c_0.$$
Hence, $U_{c_0}$ is strongly causal neighborhood of its elements.
\end{proof}

\begin{proposition}
The distance between $\tilde 1$ and $(c,w)$ with respect to the Lorentzian metric $\tilde \rho$ is equal to
$$d(\tilde 1; (c,w)) = \left\{ \begin{array}{ll}
0 & \text{if } c \geq - \tan^{-1} |w| \\
 \sin^{-1}\left( \sqrt{\tan^2 c - |w|^2} \right) & \text{if } -\frac{\pi}{2} < c < - 
\tan^{-1} |w| \\
\frac{\pi}{2} & \text{if } c = -\frac{\pi}{2} \\
\pi - \sin^{-1}\left( \sqrt{\tan^2 c - |w|^2} \right) & \text{if }  -\pi + \tan^{-1} < c 
< - \frac{\pi}{2}
\end{array}\right.$$
When $c \leq - \pi + \tan^{-1} |w|$, we know that $d(\tilde 1; (c,w)) \geq \pi$.
\end{proposition}

\begin{proof}
The identity for $c \geq -\tan^{-1} |w|$ is trivial. For $-\pi + \tan^{-1}|w| < c < 
-\tan^{-1} |w|$, the formula
follows from the fact that this space is globally hyperbolic, and from the proof of 
Proposition \ref{Prop6}.
The lower bound for the distance when $c \leq - \pi + \tan^{-1} |w|$ follows from the 
reverse triangle
inequality, and from the fact that every such element has points in its timelike past of 
distance arbitrarily close to $\pi$.

To prove that $-\pi + \tan^{-1}|w| < c < -\tan^{-1} |w|$, we need to show that for every 
pair $\tilde g_1, \tilde g_2$
in this subset, $\mathcal{J}^+(\tilde g_1, \tilde \rho) \cap \mathcal{J}^-(\tilde g_2, 
\tilde \rho)$ is compact.
It is trivial that $\mathcal{J}^+(\tilde g_1, \tilde \rho) \subset \mathcal{J}^+(\tilde 1, \tilde 
\rho)$,
and by the same reasoning $\mathcal{J}^-(\tilde g_2, \tilde \rho) \subseteq 
\mathcal{J}^-((c_0,0), \tilde \rho)$.
for some $-\pi < c_0 < 0$. Clearly,
$$\mathcal{J}^-((c_0,0), \tilde \rho) = \{(c,w) \in \wSU(1,1)| \, c \geq c_0 + \tan^{-1} |w| 
\}.$$
It follows that $\mathcal{J}^+(\tilde 1, \tilde \rho) \cap \mathcal{J}^-((c_0,0), \tilde 
\rho)$, is compact for
$- \pi < c_0 < 0$, but not when $c_0 \leq -\pi$.
\end{proof}

One would expect an analogous statement for the sub-Lorentzian metric presenting a description of
the distance function in $\Omega_0$, given the proof of Proposition \ref{maxgeo}, but
it is more difficult to prove whether it is globally hyperbolic or not.

\section{Proofs of  main results} \label{proofs}
\subsection{Proof of Proposition \ref{sRnumbergeo}}
\subsubsection{Explanation of the proof and notation}
The technique of this proof is to consider a general geodesic satisfying $\tilde g(0) = \tilde 1$
and $\tilde g(1) = (c,w)$. Each geodesic is  determined by their initial conditions
$(a_1, a_2, a_3)$, where $a_1^2 + a_2^2 \neq 0$. Our task will be, given the final point of the geodesic,
$(c,w)$, to find how many choices of initial conditions do we have.

It will be practical to use $r e^{i\theta} := a_1 + ia_2$. Notice that $\check a = r^2 - a_3^2$, is independent of $\theta$. Whenever $r \neq 0$, we will use $\beta$ to denote the value $\beta := \tfrac{a_3}{r}$. Notice that if $\beta \neq 0$, then $\check a = a_3^2 (\beta^{-2} -1)$.

We need to solve the equations
\begin{equation} \label{eqsR1}
c = - \scrT_a + a_3, \qquad
|w| = r |\scrS_{r^2 - a_3^2}| \end{equation}
\begin{equation} \label{eqsR3}
\frac{w}{|w|} = -e^{-i\theta} \frac{\scrS_{r^2 - a_3^2}}{|\scrS_{r^2 - a_3^2}|},
\qquad \text{if } w \neq 0. \end{equation}
Recall the definition of $\scrS_\alpha, \check a$ and $\scrT_a$ from Section \ref{notationsec}.

We first deal with the special case $w = 0$ right away. This implies $\scrS_{r^2 - a_3^2} = 0$ (actually, it is an if and only if), so $\sqrt{a_3^2 - r^2} = \pi k$, where $k \in \mathbb{N}$. Furthermore, $\sgn(c) = \sgn(a_3)$, $|c| = |a_3| - \pi k$ and $r = \sqrt{a_3^2 - \pi^2 k^2} = \sqrt{c(c + 2\pi k)}$,
so there are countably many choices. Also, there are no restrictions to $\theta$.

For the rest of the proof, we will assume that $|w| \neq 0$ (and hence $-\check a \neq pi^2 k^2$).
Observe that this assumption means that given values of $r$ and $a_3$, \eqref{eqsR3} determines $e^{i\theta}$.
We therefore only need to look at the number of solutions to \eqref{eqsR1}.
Since the expression for $\scrT_a$ and $\scrS_{\check a}$ vary with different choices of $a$, we need split the geodesics up into 4 different cases.

\subsubsection{If $\beta = 0$ (i.e. $a_3 = 0$)} The geodesics are contained in the horizontal
space, and there is a unique one for every point ($r = \sinh^{-1}|w|$).

\subsubsection{The case $\check a > 0, \beta \neq 0 \iff 0< |\beta| < 1$.} In order to find a solution \eqref{eqsR1}, we must solve the equations
\begin{equation} \label{betaeq11}
c = a_3 - \tan^{-1} \left(\tfrac{\beta}{\sqrt{1 - \beta^2}} \tanh \left( \tfrac{a_3 \sqrt{1 - \beta^2}}{\beta} \right) \right),
\end{equation}
\begin{equation} \label{betaeq12}
\sinh \left( \tfrac{a_3 \sqrt{1 - \beta^2} }{\beta} \right) = |w| \sqrt{1 - \beta^2}.
\end{equation}
Observe that $c - \frac{\pi}{2} < a_3 < c + \frac{\pi}{2}$.
We rewrite \eqref{betaeq12} as
$$\tanh \left( \tfrac{a_3 \sqrt{1 - \beta^2}}{\beta} \right) =
\sqrt{\tfrac{|w|^2(1 - \beta^2)}{1 + |w|^2(1 - \beta^2)}},$$
and inserting it in \eqref{betaeq11}, we obtain
$$a_3 = c + \tan^{-1}\left(\tfrac{\beta |w|}{\sqrt{1+ |w|^2(1 - \beta^2)}}\right).$$
Substituting $a_3$  back into \eqref{betaeq12}, we have that $\beta$ is a solution to the equation $B(\beta) =0$, where $B$ is defined by
\begin{equation} \label{betaeq13}
B(\beta) := \tfrac{\sqrt{1 - \beta^2}}{\beta} \left(c + \tan^{-1}\left(
\tfrac{\beta|w|}{\sqrt{1 + |w|^2(1 - \beta^2)}} \right) \right) - \sinh^{-1}(|w|\sqrt{1 - \beta^2}).
\end{equation}
This only vanishes when $0 < |c| < |w| - \tan^{-1}|w|$, and this solution is unique.

To show this, we first make the observation that the limits
$$\lim_{|\beta| \rightarrow 1-}B(\beta) = 0,\quad \lim_{|\beta| \rightarrow 0+}B(\beta) = \infty,$$
imply that a necessary condition for the existence of the solution to
\eqref{betaeq13} is that the derivative $\frac{\partial}{\partial \beta} B$ vanishes at some point. Also, if the derivative vanishes exactly once on the interval $\beta \in (0,1)$, then this is a sufficient condition for $B$ to have at most one root on this interval. The same holds for the interval $\beta \in (-1,0)$.
This derivative can be computed to be
\begin{equation} \label{derbetaeq13}
\tfrac{\partial}{\partial \beta} B(\beta) = - \tfrac{1}{\beta^2 \sqrt{1 - \beta^2}} \left(c + \tan^{-1}\left(
\tfrac{|w|\beta}{\sqrt{1 + |w|^2(1 - \beta^2)}}\right) - \tfrac{|w| \beta}{\sqrt{1 + |w|^2(1 - \beta^2)}}\right).
\end{equation}
 Assume first that $c > 0$.
Then \eqref{derbetaeq13} has no zeros for $-1 < \beta < 0$.
The derivative \eqref{derbetaeq13} vanishes at most once on the interval $0 < \beta <1$. In fact,
it vanishes exactly once only if $0 < c < |w| - \tan^{-1}|w|$, while for greater values of $c$ it does not vanish. 
With this information together with the fact that for $0< c < |w| - \tan^{-1} |w|$ it holds that $\lim_{t\to1-}\tfrac{\partial}{\partial \beta}B(\beta) <0$, we know that, in this case, \eqref{betaeq13} has one solution in  $0 < \beta <1$.

Similarly, it can be shown that \eqref{betaeq13} has exactly one solution, when
$ \tan^{-1} |w|-|w| <c<0$, and this solution $\beta$ is between -1 and 0. If $c = 0$, \eqref{derbetaeq13} never vanishes in $\beta\in (-1,1)$.

\subsubsection{The case $\check a = 0 \iff  |\beta| = 1$.} This happens when $r = |a_3|$.
It is trivial that the equations
$$c = \tan^{-1} a_3 + a_3, \qquad |w| = |a_3|,$$
have a solution, which is unique, only  if the endpoint satisfies the condition $c = \pm (|w| - \tan^{-1} |w|)$.

\subsubsection{The case $\check a > 0 \iff |\beta| > 1$:} \label{Prop2last}
For this part, we need the following lemma.

\begin{lemma} \label{boundlemma}
For $1 \leq \beta \leq \sqrt{1 + |w|^{-2}}$, define
$$\varphi_1 := \sin^{-1} \left(|w| \sqrt{\beta^2 -1}\right),
\qquad
\varphi_2 := \sin^{-1} \left(\tfrac{|w| \beta}{\sqrt{1 + |w|^2}} \right).$$
Then $(\beta + \sqrt{\beta^2 -1})(\beta \varphi_1(\beta) - \sqrt{\beta^2 -1} \varphi_2(\beta))$
monotonically increase from 0 to $\frac{\pi}{2}$ in the given interval.
\end{lemma}

\begin{proof}
Let $1<\beta < \sqrt{1+ |w|^2}.$ First, we observe that $\varphi_1$ and $\varphi_2$ have partial derivatives
$$\tfrac{\partial}{\partial \beta} \varphi_1(\beta) =
\tfrac{\beta}{\sqrt{\beta^2 -1}} \tfrac{|w|}{\sqrt{1 - |w|^2(\beta^2 -1)}}, \qquad
\tfrac{\partial}{\partial \beta} \varphi_2(\beta) =
\tfrac{|w|}{\sqrt{1 - |w|^2 (\beta^2 - 1)}}.$$
Second, the identity $$\sin^{-1} x < \sqrt{1+ |w|^2}x \quad \text{ for values } \quad 0 < x < \tfrac{|w|}{\sqrt{1+|w|^2}},$$ gives us an upper bound for $\varphi_2(\beta)-\varphi_1(\beta)$ by
\begin{align*} \varphi_2(\beta)-\varphi_1(\beta) &= \sin^{-1}\left(\tfrac{|w| (\beta - \sqrt{\beta^2 -1})\sqrt{1 -|w|^2(\beta^2-1)} }{\sqrt{1+|w|^2}}  \right) \\
& < |w| (\beta - \sqrt{\beta^2 -1}) \sqrt{1 -(\beta^2-1)|w|^2} \end{align*}
We use this to make a lower bound for the derivative and complete the proof
\begin{align*} & \tfrac{\partial}{\partial \beta} \left( (\beta + \sqrt{\beta^2 -1})
(\beta \varphi_1(\beta) - \sqrt{\beta^2 -1} \varphi_2(\beta)) \right) \\ 
= & \tfrac{\beta + \sqrt{\beta^2 -1}}{\sqrt{\beta^2 -1} \sqrt{1 - |w|^2(\beta^2 - 1)}}
\left(|w| -  \sqrt{1 - |w|^2(\beta^2 - 1)}(\beta + \sqrt{\beta^2 -1})
(\varphi_2 - \varphi_1) \right) \\
> & \tfrac{\beta + \sqrt{\beta^2 -1}}{\sqrt{\beta^2 -1} \sqrt{1 - |w|^2(\beta^2-1)}}
\Big(|w| - |w|\Big(1 -(\beta^2-1)|w|^2\Big) \Big) \\
= & \tfrac{|w|^3 (\beta + \sqrt{\beta^2 -1)}}{\sqrt{\beta^2 -1} \sqrt{1 - |w|^2(\beta^2-1)}}
(\beta^2-1) > 0.\end{align*}
\end{proof}

Equations \eqref{eqsR1} are now written in a more complicated form
\begin{equation} \label{betaeq21}
c = - \tan^{-1}\left(\tfrac{\beta}{\sqrt{\beta^2 -1}} \tan\left(\tfrac{a_3 \sqrt{\beta^2 -1}}{\beta}\right) \right)
- \pi \, \sgn(a_3) \left\lceil \tfrac{\sqrt{- \check a}}{\pi} - \tfrac{1}{2} \right\rceil + a_3,
\end{equation}
\begin{equation} \label{betaeq22}
\left|\sin\left(\tfrac{a_3 \sqrt{\beta^2 -1}}{\beta}\right) \right| = |w| \sqrt{\beta^2 -1}.
\end{equation}
An immediate observation from \eqref{betaeq22}, is that $|\beta| \leq \sqrt{1 + |w|^{-2}}$.
Let us introduce some notations. We write $\sqrt{-\check a} = \pi k + \alpha', 
\alpha' \in (- \tfrac{\pi}{2} , \tfrac{\pi}{2} ] \setminus \{ 0 \}$, and also $\sigma_1 = \sgn a_3 = \sgn \beta$,
$\sigma_2 = \sgn \alpha'$. If $k = 0$, then $\sigma_2 = 1$.
Subsequently,
$$\sin(\sqrt{-\check a}) = (-1)^{k} \sigma_2 |w| \sqrt{\beta^2 -1}, \qquad
\cos(\sqrt{-\check a}) = (-1)^{k} \sqrt{1 - |w|^2(\beta^2 - 1)}.$$
So
$$a_3 = c + \sigma_2 \tan^{-1}\left(\tfrac{|w|\beta}{\sqrt{1 - |w|^2(\beta^2 - 1)}} \right)
+ \sigma_1 \pi k = c + \sigma_2 \varphi_2(\beta) + \sigma_1 \pi k,$$
$$\alpha' = \sigma_2 \tan^{-1} \left(\tfrac{|w|\sqrt{\beta^2 -1}}{\sqrt{1 - |w|^2(\beta^2 - 1)}}\right)
= \sigma_2 \varphi_1(\beta).$$
and it follows that $\beta$ has to be a solutions to the equation
\vspace{0.3cm}
\begin{equation} \label{betaeq23}
k \pi  + \sigma_2 (|\beta| + \sqrt{\beta^2 -1})(|\beta| \varphi_1 - \sqrt{\beta^2 -1}|\varphi_2|)
=  |c| \sqrt{\beta^2 -1}(|\beta| + \sqrt{\beta^2 -1}),
\vspace{0.3cm}
\end{equation} 
where
$$\beta \in \left[-\sqrt{1 + |w|^{-2}},-1\right) \cup \left(1, \sqrt{1 + |w|^{-2}}\right] \quad \text{ and } \quad \sgn \beta = \sgn c.$$
Also $\beta$ can only have the value $\pm \sqrt{1 + |w|^{-2}}$, when $\sigma_2 = 1$. We need to count solutions for all values of $k$ and $\sigma_2 =\pm 1$. See Figure 4 for a graphical representation of the equation \eqref{betaeq23}.

From now we will just assume that $c >0$, since the considerations for $c < 0$ are totally analogous.
From Lemma \ref{boundlemma}, the left-hand  side of \eqref{betaeq23} decreases or increases from $k \pi$ to $\frac{\pi}{2} (2k + \sigma_2)$. This shows that there are no solutions for $k < 0$.
The derivative of the difference of the right and left side is given by
\begin{multline} \label{derbetaeq24}
\tfrac{\partial}{\partial \beta}\left(\left((c  + \varphi_2)\sqrt{\beta^2 -1}  - \beta \varphi_1\right) (\beta + \sqrt{\beta^2 -1})
\right) \\ 
= \tfrac{(\beta + \sqrt{\beta^2 -1})^2}{\sqrt{\beta^2 -1}}
\left(c - \tfrac{|w|(\beta - \sqrt{\beta^2 -1})}{\sqrt{1 - |w|^2(\beta^2 - 1)}}
+ \varphi_2 -\varphi_1 \right),\end{multline}
and vanishes once if $c > |w| -\tan^{-1} |w|$, and otherwise it is negative. We therefore only need to check for roots for these values of $c$, in order find roots of \eqref{betaeq23}.

If $k = 0$, then there is only a solution when $\sigma_2 =1$.
In this case, there is a unique solution if $\tfrac{|w| - \tan^{-1}|w|}{\pi (\sqrt{|w|^2 + 1} -1)}\leq \tfrac{|c|}{\pi(\sqrt{|w|^2 + 1} - 1)} \leq \tfrac{1}{2} $. The lower bound follows from \eqref{derbetaeq24}, while the upper bound is from \eqref{betaeq23}.

If $\sigma_2 = - 1$ and $k > 0$, then there is a unique root if and only if,
$k < \tfrac{|c|}{\pi(\sqrt{|w|^2 + 1} - 1)} + \frac{1}{2}$.
If $k > 0$ and $\sigma_2 = 1$, then there is clearly no solutions for
$k \geq \tfrac{|c| (\sqrt{|w|^2 + 1} +1)}{\pi |w|^2}$, and from~\eqref{derbetaeq24}, it follows that
there is only one solution for $k < \frac{|c|}{\pi(\sqrt{|w|^2 + 1} - 1)} - \frac{1}{2}$.
For the cases between these, we define the number $\hat{\beta}$ as a root of the expression
$$c - \tfrac{|w|\left(\beta - \sqrt{\beta^2 -1}\right)}{\sqrt{1 - |w|^2(\beta^2 - 1)}}
+ (\varphi_2 -\varphi_1),$$
i.e. point where the derivative in \eqref{derbetaeq24} vanishes.
Then, there is a solution  only  if
$$k \pi  + \left(\hat{\beta} + \sqrt{\hat{\beta}^2 -1}\right)\left(\hat{\beta} \varphi_1 - \sqrt{\hat{\beta}^2 -1}\varphi_2\right)
\leq  c \sqrt{\hat{\beta}^2 -1}\left(\hat{\beta} + \sqrt{\hat{\beta}^2 -1}\right).$$
If the equality is attained  in the above inequality, then $\chi_k = \tfrac{|w|\sqrt{\hat{\beta}^2 -1}}
{\sqrt{1 + |w|^2 - |w|^2 \hat{\beta}^2}}$, i.e., $\hat{\beta} = \tfrac{\sqrt{1 + |w|^{-2} + \chi_k^{-2}}}{\sqrt{1 + \chi_k^{-2}}}$.
Hence, there is only one solution when
$$|c| = \sqrt{|w|^2 + \chi_k^2 + |w|^2 \chi_k^2} - \chi_k
\tan^{-1}\left(\tfrac{\sqrt{|w|^2 + \chi_k^2 + |w|^2 \chi_k^2} - \chi_k}
{1 + \chi_k \sqrt{|w|^2 + \chi_k^2 + |w|^2 \chi_k^2}} \right),$$
and there are two solutions when $|c|$ is between this value and $\tfrac{\pi}{2}(2k +1)
(\sqrt{|w|^2 +1} - 1)$.
All the cases has now been examined, so this ends the proof.
\begin{figure}
\begin{pspicture}(-5,-0.5)(17,11)
\centering
\scalebox{0.60}{\includegraphics{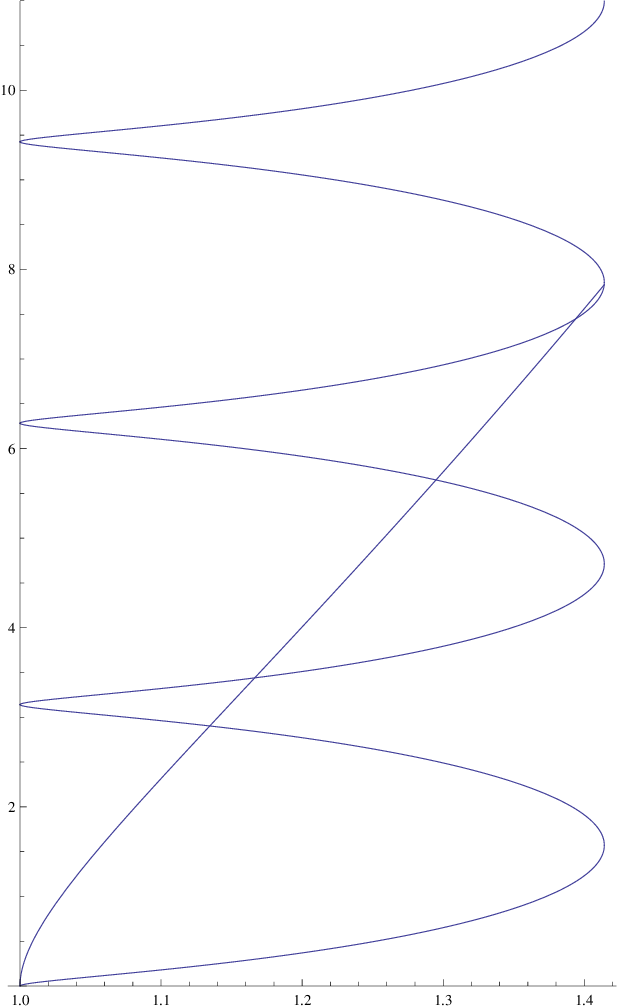}}
\rput(0.4,0.2){$\beta$}
\end{pspicture}
\caption[]{This figure essentially summarizes the results of the proof of  Proposition \ref{sRnumbergeo}.
The graph is based on the case $|w| =1$ and $|c|$ close to $\frac{5 \pi}{2}(\sqrt{1 + |w|^{-2}} -1)$, but the general idea is the for any point.
The "arches" in the graph are actually 7 different graphs, corresponding to different values of $\sigma_2$ and $k$ ($k$ is here from 0 to 3). The bottom half of the arches correspond to the case when $\sigma_2 = 1$, and the top for when $\sigma_2 = -1$. The curve that is not an arch, is the right-hand side of \eqref{betaeq23}.
If $|c|$ is such that the graph of the right-hand side of \eqref{betaeq23} goes under all the above graphs,
then the  point corresponding to this choice of parameters  is reached by a geodesic with $\check a \geq 0$,
see Corollary \ref{CorlengthsR} (c). If $|c|$ is less than
$\frac{\pi}{2}(\sqrt{|w|^{-2} +1} -1)$, but not underneath the first arch, then this corresponds to the points in
Corollary \ref{CorlengthsR} (d). For larger values, the number of geodesics is equal to the number of intersections (so, in this case, for example, there are 5 geodesics). It is also easy to see that
there can be only be either one root or no root when $\sigma_2 = -\sgn(c)$, while if $\sigma_2 = \sgn(c)$, the curve can either not intersect, be tangent (when $|c|$ is in the cut locus), or intersect twice, like the curve in our figure.}
\end{figure}
\subsection{Proof of Corollary \ref{CorlengthsR}}
Let us use the same notation as  in the proof of Proposition \ref{sRnumbergeo}.
If $\tilde g(t)$ is a geodesic that goes from $\tilde 1$ to $(c, w)$ in a unit time interval,
then $\ell(\tilde g)= r = \frac{a_3}{\beta}$. This implies (a)--(e)  rather easily from the previous proof.

To prove (f), let $g_{k,\sigma_2}$ be a geodesic corresponding to a
 choice of $k$ and $\sigma_2$, with $\beta_{k,\sigma_2}$. The length becomes
$$\ell(\tilde g_{k,\sigma_2}) = \frac{|c| + \sigma_2|\varphi_2| + k \pi}{\beta_{k,\sigma_1}}.$$
Then, from the equation $\sqrt{\beta_{k,\sigma_2}^{-2} -1} (|c| + \sigma_2 |\varphi_2| + k \pi) = \pi k + \sigma_2 \varphi_1$, we know that
$$\ell(g_{k,\sigma_2})^2 = (|c| + k\pi + \sigma_2 |\varphi_2(\beta_{k,\sigma_2})|)^2 +
(\pi k + \sigma_2 \varphi_1(\beta_{k,\sigma_2}))^2,$$
so
$$(|c| + k\pi)^2 + (\pi k)^2 < \ell(g_{k,1})^2 \leq \left(|c| + \left(k +\frac{1}{2}\right) \pi\right)^2 +
\left(\pi \left(k + \frac{1}{2}\right)\right)^2,$$
$$\left(|c| + (k- \frac{1}{2})\pi\right)^2 + \left(\pi (k- \frac{1}{2})\right)^2 \leq \ell(g_{k,-1})^2 < (|c| + k \pi)^2 +
(\pi k)^2.$$
It follows that $\tilde g_{1,-1}$ is minimal.

\subsection{Proof of Proposition \ref{maxgeo}} \label{proofmaxgeo}
\subsubsection{Explanation of the proof and notation}
Since we are working with timelike-future directed curves, we know that $a_3 >0$ and that $|a_2| < a_3$.
Let us therefore denote $a_3 = r \cosh \theta$, and $a_2 = r \sinh \theta$.
Then $\check a = a_1^2 - r^2$.
Similarly to the sub-Riemannian geodesics, we consider $\tilde g(0) = \tilde 1,$ $\tilde g(1) = (c,w)$, and look at how the final point defines the initial conditions.
We are going to use the projection of $(c,w)$ to $\SU(1,1)$. Along with the coordinates
 $z_1 = x_1 + i y_1$ and $z_2 = x_1 + i y_2$,  we shall use new coordinates given by
$$u_1 = x_1 + x_2, \quad u_2 = x_1 - x_2,
\quad v_1 = y_1 + y_2, \quad v_2 = y_1 - y_2.$$
Observe that $|z_1|^2 - |z_2|^2 = u_1 u_2 + v_1 v_2 = 1$.
Recall the definition of $\check a, \scrC_{\check a}, \scrS_{\check a}$ and $\scrT_a$ from Section \ref{notationsec}.

The projection of the geodesics to these coordinates becomes
$$\begin{array}{ll}
u_1 = e^{a_1} \left(\scrC_{\check a} - a_1 \scrS_{\check a}\right), &
u_2 = e^{-a_1}\left( \scrC_{\check a} + a_1 \scrS_{\check a} \right),\\
v_1 = - r \scrS_{\check a} e^{a_1- \theta}, &
v_2 = - r \scrS_{\check a} e^{-(a_1- \theta)}. \end{array}$$
The relationship between the projected coordinates and the original ones on the covering space is
\begin{equation} \label{formulac} \nonumber
c = - \scrT_a + \tan^{-1} \left(\frac{2\sinh a_1(u_2 v_1 - u_1 v_2) }{e^{-a_1}(u_1^2 + v_1^2 + 1) +
e^{a_1}(u_2^2 + v_2^2 +1)} \right),
\end{equation}
\begin{equation} \label{formulaw} \nonumber
w = \tfrac{1}{2}(u_1 - u_2 + i(v_1 - v_2)).
\end{equation}
Furthermore,
\begin{equation} \label{sLeq1}
2 \scrC_{a_1^2 - r^2} = u_1 e^{-a_1} + u_2 e^{a_1},
\end{equation}
\begin{equation} \label{sLeq2}
2 a_1 \scrS_{a_1^2 - r^2} = - u_1 e^{-a_1} + u_2 e^{a_1},
\end{equation}
\begin{equation} \label{sLeq3}
2 r \sinh(\theta) \scrS_{a_1^2 - r^2} = v_1 e^{-a_1} - v_2 e^{a_1},
\end{equation}
\begin{equation} \label{sLeq4}
2 r \cosh(\theta) \scrS_{a_1^2 - r^2} = - v_1 e^{-a_1} - v_2 e^{a_1}.
\end{equation}
We will use this coordinates to compute the number of geodesics. We are again forced to split into different cases.

Some of these cases we need to deal with right away. First, let us  consider the case when $a_1^2 - r^2 = - \pi^2 k^2, k \in \mathbb{N}$.
The only points that can be reached by this type of geodesics are
$(-\pi k, x_2)$. Then
$$a_1 = (-1)^k \sinh^{-1} x_2, \qquad r = \sqrt{c^2 + (\sinh^{-1}x_2)^2},$$
and there are no restrictions on $\theta$.

Otherwise, if $a_1^2 - r^2 \neq - \pi^2 k^2$, then $v_1v_2 = 1- u_1 u_1 > 0$.
It follow that
\begin{equation} \label{rthetaeq}
\theta = a_1 + \log \sqrt{\tfrac{v_2}{v_1}}, \quad
\text{and} \quad
r = \tfrac{2 \sqrt{v_1 v_2} |a_1|}{|u_2 e^{a_1} - u_1 e^{-a_1}|}  \quad
\text{if } a_1 \neq 0.
\end{equation}
Hence, when $a_1 \neq 0$, we only need to find the number of solutions to the equation \eqref{sLeq2} with respect to $a_1$. From this, the values of $r$ and $\theta$ will  be determined
by \eqref{rthetaeq}.
Remark that \eqref{rthetaeq} also give us the expression
\begin{equation} \label{checkasL} \check a = \left(\frac{a_1}{\frac{u_2}{2} e^{a_1} - \frac{u_1}{2} e^{-a_1}}\right)^2
\left(\left(\tfrac{u_2}{2} e^{a_1} + \tfrac{u_1}{2} e^{-a_1}\right)^2 - 1 \right), \quad a_1 \neq 0.\end{equation}

\subsubsection{If $a_1 = 0$:} We find  $r$ as 
$r = k + (-1)^{k+1} \sin^{-1}\left(\tfrac{(v_1 + v_2) \sqrt{v_1 v_2}}{|v_1| + |v_2|} \right)$,
where $k = \lceil \frac{c}{\pi} - \frac{1}{2} \rceil$.  There will now be exactly one geodesic, satisfying $a_1 = 0$,
when $(c,w)$ is of the form
$$c = - \pi k + \Arg(\sqrt{1 - s^2}   - \sgn(s) i \sqrt{y_2^2 + s^2}),  \quad w = i y_2,$$
where $y_2 \in \mathbb{R}, s \in (-1, 1] \setminus \{ 0 \}$ or $s \in (0,1]$ if $k =0$.

\subsubsection{If $\check a > 0$:}
A quick look at \eqref{sLeq3} and \eqref{sLeq4} yields that $v_1 < 0$ and $v_2 <0$.
Also from \eqref{sLeq1}, we know that at least one of $u_1$ and $u_2$ is positive.
We define $b = \sqrt{\left(\frac{u_2}{2} e^{a_1} + \frac{u_1}{2} e^{-a_1}\right)^2 - 1}$.
Then
$$e^{a_1} = \tfrac{\sqrt{b^2 +1} + \sigma_1 \sqrt{b^2 + s^2}}{u_2}
= \tfrac{u_1}{\sqrt{b^2 +1} - \sigma_1 \sqrt{b^2 +s^2}}, \quad \sigma_1 = \pm 1.$$
Here, we denote $s = \sigma_1\sqrt{v_1v_2}.$
Notice, that the first part of the above equation is not valid if $u_2 = 0$, and the second
part is not valid if $u_1 = 0$. The relation \eqref{sLeq2} implies $\sgn a_1 = \sigma_1$.

If $\sigma_1 = 1$, then we know that $\sgn u_2 = 1$.
Furthermore, from \eqref{sLeq2}, we obtain that $b$ is a solution to the equation
\begin{equation} \label{eqb11}
\psi_{-1}(b)-u_2=0, \qquad \mbox{where \ }
\psi_{-1}(b) := (\sqrt{b^2 +1} + \sqrt{b^2  + s^2})e^{- \frac{\sqrt{b^2 + s^2}}{b} \log(b + \sqrt{b^2 +1})},
\end{equation}
$\psi_{-1}$ increases from $(1 + s)e^{- s}$ to 1,
so a solution exists  only if $(1 + s) e^{- s} < u_2 < 1$.
Then the points, which are connected by these types of geodesics are
$$c = \Arg\left(\frac{u_2^2 + 1 - s^2}{2u_2} - i \sqrt{y_2^2 + s^2}\right),
\qquad w = \frac{1 - s^2 - u_2^2}{2u_2} + i y_2,$$
where $(1 + s) e^{-s} < u_2 < 1$.
To get some idea of the meaning the above equations,
notice that $\sgn(1- s) = \sgn u_1$, so $x_1 = \frac{1}{2}(\frac{1- s^2}{u_2} + u_2)$
is an increasing function of $u_2$ when $s > 1$. Also, $x_2 = \frac{1}{2}(\frac{1-s^2}{u_2} - u_2)$
is a decreasing function of $u_2$ when $s < 1$. 
Using this information, we know that there are only geodesic of these types connecting the points of the form
$$c = \Arg\left(x_1 - i \sqrt{s^2 + y_2^2}\right), \qquad w = x_2 + i y_2,$$
with the identity of the group, where $y_2 \in \mathbb{R}, s >0$ and $x_1, x_2$ satisfy one of the following conditions
\begin{itemize}
\item $x_2 \in (- \frac{s^2}{2}, \sinh s - s \cosh s )$,  $x_1 = \sqrt{x_2^2 + 1 - s^2}$, if $0 < s < 1$;
\item $x_1 = - x_2 \in (e^{-1}, \frac{1}{2})$, if $s = 1$;
\item $x_1 \in (\cosh s - s \sinh s, 1 - \frac{s^2}{2}), x_2 = -\sqrt{x_1^2 + s^2 -1}$, if $s > 1$.
\end{itemize} 
We can rewrite the above cases as the condition $x_2 = - \sqrt{x_1^2 + 1 - s^2}$, $
x_1 \in (\cosh s - s \sinh s, 1 - \frac{s^2}{2})$.

When $\sigma_1 = -1$, we know that $\sgn u_1 = 1$, and $b$ is a solution to the equation
$\psi_{-1}(b) = \log u_1$.
So the points reachable by these types of geodesics are
$$c = \Arg\left(x_1 - i \sqrt{y_2^2 + s^2}\right), \qquad
w = x_2 + i y_2,$$ 
where $s < 0, y_2 \in \mathbb{R}$, and $x_1, x_2$ satisfy one of the following conditions
\begin{itemize}
\item $x_2 \in (\sinh s - s \cosh s, \frac{s^2}{2} ), x_1 = \sqrt{x_2^2 + 1 - s^2}$,
if $0 > s > -1$;
\item $x_1 = x_2 \in (e^{-1}, \frac{1}{2})$ if $s = -1$;
\item $x_1 \in (\cosh s - s \sinh s, 1 - \frac{s^2}{2}), x_2 = \sqrt{x_1^2 + s^2 -1}$, if $s < -1$.
\end{itemize}

Notice that the only thing that changes with the sign of $\sigma_1$, is the sign of $x_2$.

\subsubsection{If $\check a = 0$:}
From \eqref{sLeq1} and \eqref{sLeq2}, it follows that there are geodesics only of this
type connecting the initial point with every point  of the form
$$\left\{\left. \left(\Arg\left(\cosh s - s \sinh s - i \sqrt{s^2 + y_2^2}\right),
\sinh s - s \cosh s + i y_2 \right) \right| s, y_2 \in \mathbb{R}, s \neq 0 \right\}.$$

\subsubsection{If $\check a < 0$:}
We introduce the notation $\sqrt{- \check a} = \pi k + \alpha', \alpha'
\in \left(-\tfrac{\pi}{2}, \tfrac{\pi}{2} \right] \setminus \{0 \}$. Let
$\sigma_2 = \sgn \alpha'$.
Notice that $\sgn v_1 = \sgn v_2 = (-1)^{k+1} \sigma_2$, and hence, $\sgn y_1 = (-1)^{k+1} \sigma_2$.
Now define
$$\varphi_1 := \tan^{-1}\left(\tfrac{r \cosh\theta}{\sqrt{-\alpha}} \tan \alpha'\right), \qquad
\varphi_2 := \tan^{-1} \left(\tfrac{2(u_2 v_1 - u_1 v_2)\sinh a_1 }{e^{-a_1}(u_1^2 + v_1^2 + 1) +
e^{a_1}(u_2^2 + v_2^2 +1)} \right),$$
so that $c = -\pi k  - \varphi_1 + \varphi_2$. From the fact that
$\sin c = (-1)^{k+1} \sin(\varphi_1 - \varphi_2) = y_1$, we know that
$\sgn(\varphi_1 - \varphi_2) = \sgn \varphi_1 = \sigma_2$, which yields that
$k = -\lceil \tfrac{c}{\pi} \rceil - \tfrac{1 - \sigma_2}{2}.$
Now write $c = 2n\pi + \Arg(z_1)$. Then, $2n = k + \tfrac{(1+(-1)^{k+1})\sigma_2}{2}$,
so
$$c = \left\{\begin{array}{ll} - \pi k + \Arg(- z_1), & \text{if {\it k} is odd;}\\
- k \pi + \Arg(z_1), & \text{if {\it k} is even.} \end{array} \right.$$
Given a value of $k$, the value of $c$ is completely defined
by $z_1$. We also know that if two geodesics has the same endpoint, then
the maximal difference in their value of $k$ is 1.

We need to find $a_1$ as a solution to
$$\frac{2a_1 \sqrt{1 - \left(\frac{u_2}{2} e^{a_1} + \frac{u_1}{2} e^{-a_1}\right)^2}}
{u_2 e^{a_1} - u_1 e^{-a_1}} = \sigma_2 (-1)^k \pi k + (-1)^k \sin^{-1}\left(\sqrt{1 - \left(\tfrac{u_2}{2} e^{a_1} 
+ \tfrac{u_1}{2}e^{-a_1} \right)^2}\right).$$
Let us denote $b = \sqrt{1 - \left(\frac{u_2}{2} e^{a_1} + \frac{u_1}{2} e^{-a_1}\right)^2}$. It follows that
$$e^{a_1} = \tfrac{(-1)^k (\sqrt{1 - b^2} + \sigma_1 \sqrt{v_1 v_2 - b^2})}{u_2}
= \tfrac{u_1}{(-1)^k (\sqrt{1 - b^2} - \sigma_1 \sqrt{v_1 v_2 - b^2})}, \quad \sigma_1 = \pm 1.$$
As in the previous cases, there are some obvious restrictions on the values of $u_2=0$ or $u_1=0$.
Here we used that $\sgn \left( \tfrac{u_2}{2} e^{a_1} + \tfrac{u_1}{2} e^{-a_1} \right) = (-1)^k,$
because of \eqref{sLeq1}.
From \eqref{sLeq2}, we know that $\sigma_1 = \sigma_2 \, \sgn a_1$.

If $\sigma_1 = 1$, then $\sgn u_2 = (-1)^k$, and so $b$ is a solution to
\begin{equation} \label{eqb21} \nonumber
\psi_k(b) = (-1)^k u_2, \quad
\psi_k(b) := \left(\sqrt{1- b^2} + \sqrt{s^2 - b^2}\right) e^{- \frac{\sqrt{s^2 - b^2}}{b} 
(\sigma_2 \pi k + \sin^{-1} b)},
\end{equation}
where $0 < b \leq \min\{1, |s|\}$.
Again we use the notation $s = \sigma_1\sqrt{v_1 v_2}$. The upper bound for $b$ is strict,
unless $|q_1| > 1$ and $\sigma_2 = 1$, because $\alpha'$
can be $\frac{\pi}{2}$, but not $-\frac{\pi}{2}$. Moreover, from
 \eqref{sLeq2}, $b = |s|$ implies that either $a_1 = 0$ or $\sin \sqrt{- \check a} = 0$, the case
which  already excluded. Similarly, when $\sigma_1 = -1$, so
$b$ is a solution to $\psi_k(b) = (-1)^k u_1$.

We investigate when a solution can exist.
\begin{itemize}
\item We start with the special case  $k = 0$. Then $\psi_0$ decreases from
$(1+ |s|) e^{-|s|}$ to $\sqrt{s^2 -1} e^{- \frac{\pi}{2} \sqrt{s^2 -1}}$ for
$s \geq 1$, and from $(1 + |s|) e^{-|s|}$ to $\sqrt{1- s^2}$
for $|s| <1$.
\item If $\sigma_2 = -1$, then $\psi_k$ decreases from $\infty$ to
$\sqrt{1 - s^2}$ for $|s| <1$, or to
$$\sqrt{s^2 - 1} \exp\left(\tfrac{(2k -1)\pi}{2}\sqrt{s^2 -1}\right)$$ 
otherwise.
\item If  $\sigma_2 = 1$, then $\psi_k$ increases from $0$ to $\sqrt{1 - s^2}$ for  $|s| \leq \omega_k$.
Continuing with $\sigma_2 =1$, if $|s| > \omega_k$, then $\psi_k$ increases from $0$ to
$(\sqrt{1- \omega_k^2} + \sqrt{s^2 - \omega_k^2})
\exp\left(-\sqrt{\frac{s^2 - \omega_k^2}{1 - \omega_k^2}}\right)$, and decreases to
$\sqrt{1 - s^2}$ for $|s| <1$, or to $\sqrt{s^2 - 1} \exp\left(-\frac{(2k +1)\pi}{2}\sqrt{s^2 -1}\right)$ otherwise.
\end{itemize}

Further, if $\sigma_1 = 1$, then $\sgn u_2 = (-1)^k$, and $\sgn u_1 = (-1)^k\sgn(1-s^2)$,
and the same holds when $\sigma_1 = -1$, interchanging $u_1$ and $u_2$.
Also, from $$x_1 = \frac{1}{2}\left(u_1 + \frac{1-s^2}{u_1} \right)
= \frac{1}{2} \left( \frac{1-s^2}{u_2} + u_2 \right) , \quad
x_2 = \frac{1}{2} \left( u_1 - \frac{1-s^2}{u_1} \right) = \frac{1}{2} \left( \frac{1-s^2}{u_2} - u_2 \right),$$
we know that $x_1$ is an increasing function of $u_1$ and $u_2$, when $|s| > 1$.
When $|s| < 1$, $x_2$ increases in $u_1$ and decreases in $u_2$.
As a result, we know which point is reachable by which type of geodesics.
We only show the results for $s >0$, since for $s<0$, the results are the same, only with different
sign for $x_2$.
\begin{itemize}
\item If $k = 0$, then 
\begin{itemize}
\item if $0 < s < 1$, then $x_1 = \sqrt{x_2^2 + 1 - s^2}$, $ x_2 \in
(\sinh s - s \cosh s, 0)$;
\item if $1 \leq s$, $x_2 = - \sqrt{x_1^2 + s^2 -1}$, $x_1 \in [f_+(s; 0), \cosh s - s \sinh s)$.
\end{itemize}

\item If $k > 0$ and $\sigma_2 = -1$,  then
\begin{itemize}
\item if $0 < s < 1$, $(-1)^k x_1 = \sqrt{x_2^2 + 1 - s^2}$, and $(-1)^k x_2 \in (-\infty,0)$;
\item if $s \geq 1$, $(-1)^k x_2 = - \sqrt{x_1^2 + s^2 -1}$, and $(-1)^k x_1 \in (f_-(s ; k ), \infty )$.
\end{itemize}

\item If $k > 0$ and $\sigma_2 = 1$, then there is one geodesic (list 1)
\begin{itemize}
\item if $0 < s \leq \omega_k$, $(-1)^k x_1 = \sqrt{x_2^2 + 1 - s^2}$, and $(-1)^k x_2 \in (0, \infty)$;
\item if $\omega_k < s < 1$, $(-1)^k x_1 = \sqrt{x_2^2 + 1 - s^2}$, and $(-1)^k x_2 \in [0, \infty) \newline
\cup \left\{\sqrt{1 - \omega_k^2} \sinh\left(\sqrt{\frac{s^2 - \omega_k^2}
{1 - \omega_k^2}}\right) -
\sqrt{s^2 - \omega_k^2} \cosh\left(\sqrt{\frac{s^2 - \omega_k^2}{1 - \omega_k^2}}\right) \right\}$;
\item if $s = 1$ and $(-1)^k x_1 = - (-1)^k x_2 = \sqrt{1 - \omega_k^2} e^{-1}$;
\item if $s > 1$, $(-1)^k x_2 = - \sqrt{x_1^2 + s^2 -1}$, and 
\begin{eqnarray*}(-1)^k x_1 &\in& \left(- \infty,
- \sqrt{s^2 -1} \sinh\left(\frac{(2k +1) \pi}{2} \sqrt{s^2 -1}\right)\right)\\ &\cup&
\left\{\sqrt{1 - \omega_k^2} \cosh\left(\sqrt{\frac{s^2 - \omega_k^2}
{1 - \omega_k^2}}\right) -
\sqrt{s^2 - \omega_k^2} \sinh\left(\sqrt{\frac{s^2 - \omega_k^2}{1 - \omega_k^2}}\right) \right\},
\end{eqnarray*}
\end{itemize}
and two geodesics (list 2)
\begin{itemize}
\item if $\omega_k < s < 1$, $(-1)^k x_1 = \sqrt{x_2^2 + 1 - s^2}$, and
$$(-1)^k x_2 \in \left(\sqrt{1 - \omega_k^2} \sinh\left(\sqrt{\frac{s^2 - \omega_k^2}
{1 - \omega_k^2}}\right) -
\sqrt{s^2 - \omega_k^2} \cosh\left(\sqrt{\frac{s^2 - \omega_k^2}{1 - \omega_k^2}}\right), 0\right);$$
\item if $s = 1$, and $(-1)^k x_1 = - (-1)^k x_2 \in (0, \sqrt{1 - \omega_k^2} e^{-1})$;
\item if $s > 1$, $(-1)^k x_2 = - \sqrt{x_1^2 + s^2 -1}$, and
\begin{eqnarray*}
(-1)^k x_1 \in&& \bigg[ [- \sqrt{s^2 -1} \sinh\left(\frac{(2k +1) \pi}{2} \sqrt{s^2 -1}\right), \\ && \left.
\sqrt{1 - \omega_k^2} \cosh\left(\sqrt{\frac{s^2 - \omega_k^2}
{1 - \omega_k^2}}\right) -
\sqrt{s^2 - \omega_k^2} \sinh\left(\sqrt{\frac{s^2 - \omega_k^2}{1 - \omega_k^2}}\right)\right).
\end{eqnarray*}
\end{itemize}
\end{itemize}
The above list can be rewritten so the interval is determined for $x_1$ only.
Note that the term "one geodesic" and "two geodesics" is a bit misleading,
since, if we include both signs of $s$, list 1 and list 2 are not disjoint.

\subsubsection{Summary}
We start by defining some sets to sum up the information obtained so far:
\begin{itemize}
\item $\tilde A_0$ consist of all points
$$c = \Arg\left(x_1 - i \sqrt{s^2 + y_2^2}\right), \quad w = \pm \sqrt{x_1^2 + s^2 -1} + i y_2,$$
where $s, y_2 \in \mathbb{R}$, $s > 0$, $x_1 \in (f_+(s;0) , 1 - \frac{s^2}{2})$, when $0<  s \leq 1$,
and $x_1 \in [f_+(s;0) , 1 - \frac{s^2}{2})$, when $s > 1$.
They are the points that are reachable by geodesics with $\alpha_a \geq 0$,
and with $\alpha_a < 0$, $k = 0$.
\item $A_k^0, k \in \mathbb{N}_0$, consist of the points
$$c = - \pi k + \Arg(\sqrt{1 - s^2}   -  i \sqrt{y_2^2 + s^2}),  \quad w = i y_2,$$
where $y_2 \in \mathbb{R}, s \in (0,1]$, and  define
$A_{-k}^0, k \in \mathbb{N}$, by.
$$c = - \pi k + \Arg(\sqrt{1 - s^2}   +  i \sqrt{y_2^2 + s^2}),  \quad w = i y_2,$$
where $y_2 \in \mathbb{R}, s \in (0, 1)$.
They are the points reachable by geodesics with $a_1 = 0$.
\item $A_k^-, k \in \mathbb{N}$, consist of all points
$$c = - (k-1) \pi + \Arg\left(x_1 - i \sqrt{s^2 + y_2^2}\right), \quad
w = \pm \sqrt{x_1^2 + s^2 -1} + i y_2,$$
where $s, y_2 \in \mathbb{R}$, $s > 0$ and $x_1 \in (-\infty, -f_-(s; k))$.
They are the points reachable by geodesic with a certain choice of $k$, and with $\sigma_2 = -1$.
\item $A_k^+, k \in \mathbb{N}$, consist of all points
$$c = - k \pi + \Arg\left(x_1 - i \sqrt{s^2 + y_2^2}\right), \quad
w = \pm \sqrt{x_1^2 + s^2 -1} + i y_2,$$
where $s, y_2 \in \mathbb{R}$, $s > 0, s \neq 1$, and
\begin{itemize}
\item $x_1 \in (- \infty, f_+(s;k)) \cup \{F(s,\omega_k) \}$, when $s > 1$;
\item $x_1 = F(1,\omega_k)$ if $s = 1$;
\item $x_1 \in [f_+(s;k), \infty) \setminus \{F(s, \omega_k \}$, when $\omega_k< s < 1$;
\item $x_1 \in (f_+(s;k), \infty)$, when $0< s \leq \omega_k$.
\end{itemize}
They are the points from list 1, for $\sigma_2 = 1$ and for a choice of $k$
(except for the points $x_1 = F(s, \omega_k)$, for $\omega_k < s < 1$,  that
appear twice and which are therefore included in $\tilde B_k$).
\item $\tilde B_k, k \in \mathbb{N}$, consist of all points
$$c = - k \pi + \Arg\left(x_1 - i \sqrt{s^2 + y_2^2}\right), \quad
w = \pm \sqrt{x_1^2 + s^2 -1} + i y_2,$$
where $s, y_2 \in \mathbb{R}$, $s > \omega_k$ and
\begin{itemize}
\item  $x_1 \in \left(f_+(s; k) , F(s,\omega_k)\right]$ when $\omega_k < s < 1$;
\item $x_1 \in \left(0 , F(1,\omega_k)\right)$ when $s = 1$;
\item $x_1 \in [f_+(s;k) , F(s, \omega_k))$ when $s >1$.
\end{itemize}
They are the points from list 2 for a certain choice of $k$ and with $\sigma_2 = 1$.
\end{itemize}
Notice that $A_0^0$ is on the boundary of $\tilde A_0$, $A^0_{-k}$ is on the boundary of $A_k^-$
and $A^0_k$ is on the boundary of $A^+_k$ for $k >0$. 
The result follows  merely by comparing the sets and counting number of the sets in which
a value of $x_1$ appears. Points in $\tilde B_k$ are counted twice. We use
$\tilde A_0$, $A_0^0$, $ A^0_{-1}$, $ A_1^-$, and $A_{-k}^0$ to define $\Omega_0$, and use $A_k^+$, $A_k^0$,
 $A_{k+1}^-$, $ A_{-k-1}^0$, and $\tilde B_k$ to define $A_k$, $ B_k$, and $C_k$.

This ends the proof.


\begin{thebibliography}{99}

\bibitem{AExp}
A.~Agrachev, {\it Exponential mapping for contact sub-Riemannian structures},
Journal of Dynamical and Control Systems, {\bf  2}, (1996), no. 3, 321--358. 

\bibitem{AS}
A.~Agrachev, Yu.~Sachkov, {\it Control theory from the geometric viewpoint}, 
Encyclopaedia of Math. Sci., 87. Control Theory and Optimization, II. Springer-Verlag, Berlin, 2004,~412 pp.

\bibitem{BelsR}
A.~Bella{\"\i}che, J.-J.~Risler, {\it Sub-Riemannian Geometry}, Birkh\"auser Verlag, Basel, 1996, 393~pp.

\bibitem{BoRo}
U.~Boscain, F.~Rossi, {\it Invariant Carnot-Caratheodory metric on $S^3$, $SO(3)$, $SL(2)$ and lens spaces},
Siam J. Control Optim. {\bf 47} (2008), no. 4,  1851--1878.

\bibitem{CalinChang}
O.~Calin, D.-C.~Chang, {\it Sub-Riemannian geometry. General theory and examples}, Cambridge Univ. Press, 2009, 384~pp.

\bibitem{CChM}
O.~Calin, D.-C.~Chang, I.~Markina, {\it Sub-Riemannian geometry
of the sphere $S^3$},  Canadian J. Math. {\bf 61} (2009), no. 4, 821--839.

\bibitem{Car}
S.~Carlip, {\it Conformal field theory, $(2+1)$-dimensional gravity and the BTZ black hole},  Classical Quantum Gravity  {\bf 22}  (2005),  no. 12, R85--R123. 

\bibitem{ChangMarkVas}
D.-C.~Chang, I.~Markina, A.~Vasil'ev, {\it Sub-Lorentzian geometry on anti-de Sitter space},  J. Math. Pures Appl. {\bf 90} (2008), no.1, 82--110.

\bibitem{ChangMarkVas2}
D.-C.~Chang, I.~Markina, A.~Vasil'ev, {\it Hopf fibration: geodesics and distances}, J. Geom. Phys., {\bf 61} (2011), 986--1000.


\bibitem{Chow}
W.~L.~Chow. {\it \"Uber Systeme von linearen partiellen Differentialgleichungen erster Ordnung}, Math. Ann.,
{\bf 117} (1939), 98--105.

\bibitem{Groch1}
M.~Grochowski, {\it Geodesics in the sub-Lorentzian geometry},  Bull. Polish Acad. Sci. Math.  {\bf 50}  (2002),  no. 2, 161--178. 

\bibitem{Groch2}
M.~Grochowski, {\it On the Heisenberg sub-Lorentzian metric on $\mathbb R^3$},  Geometric Singularity Theory, Banach Center Publ., Polish Acad. Sci. {\bf 65} (2004), 57--65.

\bibitem{Groch3}
M.~Grochowski, {\it Reachable sets for the Heisenberg sub-Lorentzian structure $\mathbb R^3$.
An estimate for the distance function},  J. Dynamical and Control Systems, {\bf 12} (2006), no. 2, 145--160.



\bibitem{Jur1}
V.~Jurdjevic, {\it Geometric control theory}, Cambridge Studies in Adv. Math., 52. Cambridge Univ. Press, Cambridge, 1997,~492~pp.


\bibitem{Korolko}
A.~Korolko, I.~Markina, {\it Nonholonomic Lorentzian geometry on some $\mathbb H$-type groups},  J. Geom. Anal. {\bf 19 }(2009), 864--889.

\bibitem{LiuSussman}
W.~Liu, H.~J.~Sussman, {\it Shortest paths for sub-Riemannian metrics on rank-two distributions},
Mem. Amer. Math. Soc.  {\bf 118}  (1995),  no. 564,~104 pp.

\bibitem{Mon}
R.~Montgomery, {\it A tour of subriemannian geometries, their geodesics and applications}, Mathematical Surveys and Monographs, {\bf 91}. American Mathematical Society, Providence, RI, 2002, 259~pp. 


\bibitem{Rashevsky}
P.~K.~Rashevski{\u\i}, {\it About connecting two points of complete nonholonomic space by admissible curve}, Uch. Zapiski Ped. Inst. K.~Liebknecht {\bf 2} (1938), 83--94.


\bibitem{Str1}
R.~S.~Strichartz, {\it Sub-Riemannian geometry},  J. Differential Geom. {\bf 24}  (1986),  no. 2, 221--263.

\bibitem{Str2}
R.~S.~Strichartz, {\it Corrections to: "Sub-Riemannian geometry" J. Differential Geom. {\bf 24} (1986), no. 2, 221--263;} J. Differential Geom.  {bf 30}  (1989),  no. 2, 595--596. 


\bibitem{VershikGersh}
A.~M.~Vershik, V.~Ya.~Gershkovich, {\it Geodesic flows on $\SL(2,\mathbb R)$ with nonholonomic restrictions}, Zap. Nauchn. Semin. LOMI, {\bf 155} (1986), 7--17.

\bibitem{Witten}
E.~Witten, {\it String theory and black holes}, Phys. Rev. D {\bf 44} (1991), no. 2, 314--324.

\end{thebibliography}
\end{document}